\documentclass{cmip-l}
\usepackage{amsfonts,amssymb,amsmath}
\usepackage[french]{babel}
\def\up#1{\raise 1ex\hbox{\small #1}}
\newcommand{\cosa}{\!\!\!\!\!\!}

\newcommand{\remark}{\subsection{Remarque.}}
\newcommand{\theorem}[2]{\subsection{Th\'eor\`eme #1.} \itshape{#2}}

\newcommand{\lemma}[1]{\subsection{Lemme.} \itshape{#1}}

\newcommand{\demo}{\\ \\ \textsl{D\'emonstration.} }

\usepackage[all]{xy}
\usepackage[pdfpagemode=FullScreen,bookmarks=true]{hyperref}
\title{LE TH\'EOR\`EME DE P\'ERIODICIT\'E EN $K$-TH\'EORIE HERMITIENNE}
\author{Max Karoubi}
\begin{document}
\maketitle

La p\'eriodicit\'e de Bott joue un r\^ole primordial en $K$-th\'eorie topologique. Elle est d'ailleurs li\'ee
intimement au th\'eor\`eme d'Atiyah-Singer et plus g\'en\'eralement \`a la g\'eom\'etrie non commutative.
Dans deux articles pr\'ec\'edents \cite{K1} et \cite{K2}, nous avons d\'emontr\'e l'analogue de ce th\'eor\`eme en
$K$-th\'eorie hermitienne pour des anneaux \underline{discrets} avec (anti)involution $a {\mapsto}\overline{a}$, sous l'hypoth\`ese
qu'il existe un \'el\'ement $\lambda$  du centre de $A$ tel que $\lambda  +\overline{\lambda}  =1$ (on dit alors que 1 est scind\'e dans
$A$). Si l'anneau est commutatif et muni de l'involution triviale, ceci introduit l'hypoth\`ese que
2 est inversible dans $A$.\\
Si cette derni\`ere hypoth\`ese est anodine pour les alg\`ebres de Banach, il n'en est pas de m\^eme
pour des anneaux importants comme l'anneau de groupe $\mathbb{Z}{\pi}$, o\`u ${\pi}$ est un groupe discret. Une
difficult\'e rencontr\'ee pour l'\'etude de ce type d'anneau est la divergence entre les notions de
forme quadratique et de forme hermitienne. Dans cet article, nous d\'eveloppons une th\'eorie qui d\'epasse cette dichotomie et qui est d\'ej\`a pr\'esente dans le travail fondamental de Ranicki \cite{R}. Gr\^ace \`a cette th\'eorie le th\'eor\`eme de p\'eriodicit\'e peut \^etre d\'emontr\'e pour tout anneau. Nous montrons par exemple que les groupes de Witt sup\'erieurs d'un corps fini de caract\'eristique 2 sont tous isomorphes \`a $\mathbb{Z}/2$ (exemple 5.14).
\\
Les m\'ethodes de cet article sont beaucoup inspir\'ees de celles de \cite{K1} et \cite{K2} que nous
adaptons \`a notre propos, ce qui nous permet d'\^etre relativement bref pour certaines
d\'emonstrations. Un autre ingr\'edient essentiel est un cup-produit entre formes quadratiques
d\'efini par Clauwens \cite{C} . Celui-ci permet de d\'efinir le morphisme de p\'eriodicit\'e dans le cas g\'en\'eral. 
L'article de Clauwens ayant \'et\'e \'ecrit dans un contexte diff\'erent, nous reprenons dans
un appendice les lemmes essentiels dont nous avons besoin pour nos d\'emonstrations.\\
R\'esumons bri\`evement les diff\'erentes parties de cet article

\begin{enumerate}
\item \textbf{Description de diff\'erents types de formes hermitiennes.}
Apr\`es des rappels sur les d\'efinitions classiques utilis\'ees, nous introduisons un nouveau type de
groupe orthogonal, dit ``\'elargi'' : cf. 1.6/7. Si 1 est scind\'e dans $A$, celui-ci co\"{i}ncide avec le
groupe orthogonal sur l'anneau des nombres duaux associ\'e \`a $A$, soit $A[e]/e^2$, not\'e simplement
$A(e)$ dans la suite de l'article.
\item \textbf{Les groupes de Grothendieck et Bass en $K-$th\'eorie hermitienne.}
Nous montrons comment les th\'eor\`emes principaux en $K-$th\'eorie hermitienne restent valables
dans le cas ``\'elargi''. Nous pr\'ecisons aussi les notations utilis\'ees, en suivant partiellement la
terminologie du livre de Bak \cite{B}. Par exemple, la notation ``$L$'' , utilis\'ee en \cite{K1} et \cite{K2}, est
abandonn\'ee et remplac\'ee par la notation ``$\mathcal{K\mathcal{Q}}$'' , pour \'eviter toute ambig\"{u}it\'e avec les groupes de
chirurgie.
\item \textbf{Les groupes ${}_{\varepsilon}\bf K\mathcal{Q}_{n}(A)$ pour $\bf n > 0$ et $\bf n < 0$.}
Les d\'efinitions essentielles sont contenues dans ce paragraphe, en utilisant des id\'ees bien
connues en $K-$th\'eorie alg\'ebrique. Le th\'eor\`eme 3.2 permet de comparer les th\'eories ``max'' et
``min'', suivant la terminologie de Bak. Nous montrons aussi comment les techniques de Quillen
se transcrivent dans notre situation en une description plus g\'eom\'etrique des \'el\'ements de 
${}_{\varepsilon}K\mathcal{Q}_{n}(A)$.
\item \textbf{Cup-produits en $K-$th\'eorie hermitienne. Le cup-produit de Clauwens.}
Le cup-produit en $K-$th\'eorie hermitienne est d\'efini \`a l'aide de sa description en termes de fibr\'es
plats. Un cup-produit plus subtil, d\^{u} essentiellement \`a Clauwens, est d\'efini en 4.3 (cf. aussi
l'appendice). Nous montrons comment tous ces produits sont reli\'es entre eux dans le th\'eor\`eme 4.7.
\item \textbf{Le th\'eor\`eme fondamental de la $K-$th\'eorie hermitienne pour des anneaux arbitraires.}
Dans ce paragraphe, nous g\'en\'eralisons les r\'esultats principaux de \cite{K1} et \cite{K2} (cf. le th\'eor\`eme
5.2 et la remarque 5.11). La relation avec les groupes de Witt est faite dans le th\'eor\`eme 5.10.
\item \textbf{Les groupes de Witt stabilis\'es.}
En utilisant les r\'esultats pr\'ec\'edents, nous introduisons une th\'eorie nouvelle de groupes de Witt
``stabilis\'es'' g\'en\'eralisant ceux d\'efinis en \cite{K4}. Ses propri\'et\'es fondamentales sont d\'ecrites en
6.1. Une g\'en\'eralisation dans le cadre des sch\'emas a \'et\'e propos\'ee par M. Schlichting \cite{S} en supposant 
2 inversible.
\item\textbf{Appendice.} Les lemmes de Clauwens.
\end{enumerate}

\textbf{Remerciements.} Ce travail a \'et\'e essentiellement accompli pendant le programme th\'ematique
sur la th\'eorie de l'homotopie en 2007, organis\'e au Fields Institute \`a Toronto. Je remercie
\'egalement A. Ranicki pour avoir attir\'e mon attention sur l'article de Clauwens \cite{C}, J. Berrick
pour la d\'emonstration du lemme 4 en appendice, plus simple que le lemme original de
Clauwens, ainsi que M. Schlichting pour des commentaires pertinents apr\`es une premi\`ere
version de ce texte.

\section{Description des diff\'erents types de formes hermitiennes et quadratiques.}
\subsection{} Soit $A$ un anneau muni d'une (anti)involution $a{\mapsto}\overline{a}$ 
(on dit alors que $A$ est un anneau
hermitien) et soit\footnote{On pourrait choisir plus g\'en\'eralement un \'el\'ement ${\varepsilon}$ du centre de 
$A$ tel que ${\varepsilon}\overline{\varepsilon}=1$. Cependant, on se ram\`ene \`a
ce cas en rempla\c{c}ant $A$ par $M_2(A)$, l'alg\`ebre des matrices $2\times 2$ \`a coefficients dans $A$, 
munie d'une involution ad\'equate (cf. 1.10).} ${\varepsilon}={\pm}1$ . Nous d\'esignons par $\mathcal{P}(A)$ 
la cat\'egorie des $A-$modules (\`a droite) qui sont projectifs de type fini (les morphismes \'etant restreints 
aux isomorphismes). Si $E$ est un objet de $\mathcal{P}(A)$, son dual $E^*$ est le groupe form\'e des 
applications additives $f : E \to A$ telles que
$f(x{\lambda} ) = \overline{\lambda} f(x)$, o\`{u} ${\lambda} \in  A$ et $x \in E$. C'est en fait un objet de 
$\mathcal{P}(A)$, la structure de $A$-module \`a droite \'etant d\'efinie par la formule 
$(f.{\lambda}) (x) = f(x){\lambda}$ . Le module $E$ et son bidual $E^{**}$ sont
isomorphes canoniquement gr\^{a}ce \`a la correspondance $x {\mapsto} ( f {\mapsto} \overline{f (x)} )$. 
Nous identifierons $E$ \`a $E^{**}$ par cet isomorphisme. Par ailleurs, si $f : E \to F$ est 
un morphisme dans $\mathcal{P}(A)$, son transpos\'e ${}^{t}f : F^* \to E^*$ est d\'efini par la 
formule classique ${}^{t}(f)(g) = g.f$ et on a ${}^t({}^{t}f) =f$, compte tenu des isomorphismes canoniques 
entre les modules $E$, $F$ et leurs biduaux respectifs.

\subsection{} Nous d\'efinissons une forme ${\varepsilon}$-hermitienne sur $E$ comme un morphisme
${\phi}:E\to E^*$ tel que ${}^t{\phi}={\varepsilon}{\phi}$, o\`{u} ${}^t{\phi}:(E^*)^* \cong E \to E^*$. 
La forme ${\phi}$ est dite ``non d\'eg\'en\'er\'ee'' si c'est un isomorphisme.
Il convient de remarquer que la donn\'ee de ${\phi}$ \'equivaut \`a celle d'une application $\mathbb{Z}-$bilin\'eaire
$${\chi}:E\times E\to A$$
telle que ${\chi} (x{\lambda} , y {\mu} )=\overline{\lambda}{\chi}(x,y){\mu}$ si ${\lambda}$ et ${\mu} \in  A$, 
$x$ et $y \in E$. La correspondance est donn\'ee par la formule classique suivante
$${\chi} (x, y) = {\phi} (y)(x)$$
La condition de ${\varepsilon}$-symetrie (${}^t{\phi}={\varepsilon}{\phi}$) se traduit par l'identit\'e
$${\chi} (y, x) = {\varepsilon} \overline{{\chi} (x,y)}$$
Dans cet article, les formes hermitiennes ${\phi}$  qui nous int\'eressent sont paires : elles s'\'ecrivent
sous la forme
$${\phi}  = {\phi}_0 + {\varepsilon}{}^t{\phi}_0$$
Il convient de noter que ${\phi}_0$ n'est pas d\'etermin\'e par cette formule. Si ${\phi}_1$ est un autre choix et si
pose ${\gamma} = {\phi}_0 - {\phi}_1$, on a $^t{\gamma} = - {\varepsilon}{\gamma}$.
\subsection{} Les formes hermitiennes paires sont les objets d'une cat\'egorie 
not\'ee\footnote{ En suivant la terminologie de Bak \cite{B}.} ${}_{\varepsilon}\mathcal{Q}^{\hbox{\tiny{max}}}(A)$, d\'efinie
de la mani\`ere suivante : un morphisme
$$(E, {\phi} ) \to (F, {\psi} )$$
est un isomorphisme $f$ entre les $A$-modules sous-jacents tel que le diagramme suivant commute
$$\xymatrix{E\ar[d]_{\phi} \ar[r]^f& F\ar[d]^{\psi} \\
E^*&\ar[l]^{{}^t f} F^*}$$
\subsection{} De mani\`ere parall\`ele, en suivant Tits \cite{T} et Wall \cite{W}, on d\'efinit une forme 
${\varepsilon}$-quadratique non d\'eg\'en\'er\'ee sur $E$ comme une classe de morphismes
$${\phi}_0 : E \to E^*$$
tels que ${\phi}_0 + {\varepsilon}{}^t{\phi}_0 = {\phi}$  soit une forme hermitienne non d\'eg\'en\'er\'ee. 
Plus pr\'ecis\'ement, la classe de ${\phi}_0$ est d\'efinie modulo l'addition par un morphisme du type 
${}^t{\gamma}-{\varepsilon}{\gamma}$. Les formes ${\varepsilon}$-quadratiques sont aussi les objets d'une 
cat\'egorie not\'ee\footnote{cf. la note pr\'ec\'edente.} 
${}_{\varepsilon}\mathcal{Q}^{\hbox{\tiny{min}}}(A)$. 
Un morphisme
$$(E, {\phi}_0) \to (F, {\psi}_0)$$
est un isomorphisme $f$ entre les $A$-modules sous-jacents tel qu'il existe ${\gamma}$, morphisme de $E$ dans
$E^*$, v\'erifiant l'identit\'e 
$${}^{t}f.{\psi}_0.f = {\phi}_0 + {\gamma} - {\varepsilon} {}^t{\gamma}\qquad(S)$$
\subsection{Remarques.} Si $A$ est un corps muni de l'involution triviale, il est facile de voir que la
cat\'egorie des 1-formes quadratiques est \'equivalente \`a la cat\'egorie usuelle : il suffit de poser
$$q(x) = {\phi}_0(x)(x)$$
Cette remarque justifie la d\'efinition abstraite introduite dans 1.3.\\ \\
Par ailleurs, si 1 est scind\'e dans $A$ (cf. l'introduction), la cat\'egorie des modules 
${\varepsilon}$-hermitiens est \'equivalente \`a celle des modules ${\varepsilon}$-quadratiques : 
avec les d\'efinitions ci-dessus il suffit de poser
${\gamma} = {\lambda} ({}^{t}f.{\psi}_0.f - {\phi}_0)$. Ce cas se pr\'esente notamment si 2 est inversible dans $A$.

\subsection{}  Nous allons maintenant introduire une troisi\`eme cat\'egorie qui jouera un r\^ole important dans
notre travail et qui sera not\'ee $\mathcal{Q}^{\hbox{\tiny{\'el}}}(A)$ (``\'el'' pour ``\'elargi'' ; cf. la fin de 1.7). Les objets sont
quasiment les m\^emes que ceux de la cat\'egorie ${}_{\varepsilon}\mathcal{Q}^{\hbox{\tiny{min}}}(A)$ pr\'ec\'edente, 
sauf que l'on consid\`ere les ${\phi}_0$ comme donn\'es dans la structure (on ne consid\`ere pas 
seulement les \underline{classes} de tels ${\phi}_0$). Un morphisme de $(E, {\phi}_0)$ vers $(F, {\psi}_0)$ est d\'efini 
par un couple $(f, {\gamma})$, tel que l'identit\'e (S) ci-dessus soit satisfaite.
La loi de composition des morphismes s'explicite ainsi
$$(f, {\gamma} ).(g, {\zeta}) = (f.g, {\zeta} + {}^tg.{\gamma}.g)\qquad (C)$$
ce qui est coh\'erent avec l'identit\'e $(S)$.

\subsection{}  Il est instructif de d\'ecrire plus pr\'ecis\'ement le groupe des automorphismes d'un objet dans
chacune des trois cat\'egories. Si $(E, {\phi} )$ est un objet de ${}_{\varepsilon}\mathcal{Q}^{\hbox{\tiny{max}}}(A)$, 
\textbf{le groupe unitaire} ${}_{\varepsilon}O^{\hbox{\tiny{max}}}(E, {\phi} )$ est d\'efini par des isomorphismes $f : E \to E$ 
tels que
$${}^{t}f.{\phi} .f = {\phi}$$
Si on note $f^* = {\phi}^{-1}. {}^{t}f.{\phi}$  l'op\'erateur adjoint de $f$, il revient au m\^eme d'\'ecrire
$$f^*. f = Id_E\hbox{ (ou $f. f^* = Id_E$)}$$
Le \textbf{groupe orthogonal} ${}_{\varepsilon}O^{\hbox{\tiny{min}}}(E, {\phi}_0)$ est d\'efini par des 
isomorphismes $f : E \to E$ tels qu'il existe ${\gamma}$, morphisme de $E$ dans $E^*$, v\'erifiant l'identit\'e
$${}^{t}f.{\phi}_0.f = {\phi}_0 + {\gamma} - {\varepsilon}{}^t{\gamma}\qquad (E)$$
Il est clair que ${}_{\varepsilon}O^{\hbox{\tiny{min}}}(E, {\phi}_0)$ est un sous-groupe de 
${}_{\varepsilon}O^{\hbox{\tiny{max}}}(E, {\phi})$ pour ${\phi}  = {\phi}_0 + {\varepsilon}{}^t{\phi}_0$, la forme
hermitienne associ\'ee \`a ${\phi}_0$. Il est facile de voir que les groupes 
${}_{\varepsilon}O^{\hbox{\tiny{max}}}(E,{\phi})$ et ${}_{\varepsilon}O^{\hbox{\tiny{min}}}(E, {\phi}_0)$ 
co\"{i}ncident si 1 est scind\'e dans $A$.\\ \\
Finalement, le \textbf{groupe orthogonal \'elargi} ${}_{\varepsilon}O^{\hbox{\tiny{\'el}}} (E, {\phi}_0)$ est 
d\'efini par des couples $(f, {\gamma})$ v\'erifiant l'identit\'e $(E)$ ci-dessus. La loi de composition 
est donn\'ee par l'identit\'e $(C)$ \'ecrite aussi plus haut. On a un \'epimorphisme 
$${}_{\varepsilon}O^{\hbox{\tiny{\'el}}} (E, {\phi}_0) \to {}_{\varepsilon}O^{\hbox{\tiny{min}}}(E, {\phi}_0)$$
dont le noyau est \'egal au groupe ab\'elien ${}_{\varepsilon}S(E)$ form\'e des morphismes 
${\gamma} : E \to E^*$ tels que ${}^t{\gamma} ={\varepsilon} {\gamma}$. Nous 
obtenons ainsi une extension de groupes non triviale en g\'en\'eral
$$\xymatrix{1\ar[r]&{}_{\varepsilon}S(E) \ar[r]& {}_{\varepsilon}O^{\hbox{\tiny{\'el}}} (E, {\phi}_0) \ar[r]&
{}_{\varepsilon}O^{\hbox{\tiny{min}}}(E, {\phi}_0) \ar[r]& 1}$$
Celle-ci justifie la terminologie adopt\'ee de ``groupe orthogonal \'elargi''.\\ \\
Pour les calculs, il est commode d'identifier $E$ \`a son dual par la forme hermitienne ${\phi}$ associ\'ee \`a
${\phi}_0$. Le morphisme ${\gamma}$ est alors remplac\'e par un endomorphisme $u = {\phi}^{-1}{\gamma}$ de $E$. On peut de
m\^eme remplacer ${\phi}_0$ par ${\psi}={\phi}^{-1}{\phi}_0$. On a alors 
${\psi}^* ={\phi}^{-1}.{}^t{\psi}.{\phi}={\phi}^{-1}({}^t {\phi}_0. \overline{\varepsilon} {\phi}^{-1}.{\phi})={\phi}^{-1}({\phi} - {\phi}_0) =1-{\psi}$. 
La relation $(E)$ ci-dessus s'\'ecrit alors $f^*.{\psi}.f={\psi}+u-u^*$ ou encore $f^{-1}.{\psi}.f =
{\psi} + u - u^*$ puisque $f$ est unitaire.\\ \\
Gr\^{a}ce \`a cette traduction, la loi de composition dans ${}_{\varepsilon}O^{\hbox{\tiny{\'el}}} (E, {\phi}_0)$ 
s'\'ecrit simplement
$$(f, u).(g, v) = (f.g, v + g^{*}.u.g) = (f.g, v + g^{-1}.u.g)$$
Le noyau ${}_{\varepsilon}S(E)$ de l'homomorphisme surjectif 
${}_{\varepsilon}O^{\hbox{\tiny{\'el}}}(E,{\phi}_0)\to{}_{\varepsilon}O^{\hbox{\tiny{min}}}(E,{\phi}_0)$ 
s'identifie \`a l'ensemble des morphismes auto-adjoints de $E$, not\'e simplement $S(E)$. L'extension pr\'ec\'edente
s'\'ecrit alors de mani\`ere \'equivalente
$$\xymatrix{1 \to S(E) \ar[r]& {}_{\varepsilon}O^{\hbox{\tiny{\'el}}} (E) \ar[r]& 
{}_{\varepsilon}O^{\hbox{\tiny{min}}}(E) \ar[r]& 1}$$
Dans cette extension, le groupe ${}_{\varepsilon}O^{\hbox{\tiny{min}}}(E)$ op\`ere \`a droite sur $S(E)$ 
par la formule suivante :
$$(u, g) {\mapsto} g^{-1}.u.g$$

\subsection{}  Si 1 est scind\'e dans $A$, on peut d\'efinir une section $s$ de cette extension en posant
$$s(g) = (g, {\lambda}(g^{*}.{\psi} .g- {\psi} )) = (g, {\lambda}(g^{-1}.{\psi} .g - {\psi} ))$$
Il en r\'esulte que le groupe orthogonal \'elargi s'identifie au produit semi-direct du groupe
orthogonal ${}_{\varepsilon}O(E)$ par le groupe additif $S(E)$, gr\^{a}ce \`a l'action d\'efinie ci-dessus.
Une autre fa\c{c}on de voir les choses est d'introduire l'anneau des nombres duaux $A(e)$ avec 
$e^2 =0$ et $\overline{e}= - e$ puis d'\'etendre les scalaires \`a $A(e)$. Nous savons d\'ej\`a que le groupe orthogonal
${}_{\varepsilon}O^{\hbox{\tiny{min}}}(E)$ s'identifie au groupe unitaire 
${}_{\varepsilon}O^{\hbox{\tiny{max}}}(E)$. Par ailleurs l'\'epimorphisme
$$O^{\hbox{\tiny{max}}}(E(e)) \to O^{\hbox{\tiny{max}}}(E)$$
a comme noyau l'ensemble des matrices unitaires du type $1 + u e$, c'est-\`a-dire v\'erifiant l'identit\'e
$(1 + ue)(1 - u^*e) = 1 + (u - u^*)e = 1$, soit $u = u^*$. Le groupe unitaire op\`ere sur ce noyau par
l'action \`a droite d\'efinie par la m\^eme action : $(u, g) {\mapsto} g^{-1}u g$. Il en r\'esulte que le groupe
orthogonal \'elargi $O^{\hbox{\tiny{\'el}}}(E)$ s'identifie \`a $O^{\hbox{\tiny{max}}}(E(e))$ 
en tant que produit semi-direct.

\subsection{} Consid\'erons le cas particulier o\`{u} $A = B \times B^{op}$, $B^{op}$ \'etant l'anneau oppos\'e \`a 
$B$, l'involution permutant les facteurs du produit. Si nous posons ${\lambda} = (1, 0)$, on a 
${\lambda} + \overline{\lambda} = 1$, ce qui montre
que 1 est scind\'e dans $A$. Il est facile de voir que la donn\'ee d'un $A$-module hermitien \'equivaut \`a
celle d'un $B$-module. Les cat\'egories ${}_{\varepsilon}\mathcal{Q}^{\hbox{\tiny{max}}}(A)$ et 
${}_{\varepsilon}\mathcal{Q}^{\hbox{\tiny{min}}}(A)$ sont donc toutes les deux
\'equivalentes \`a la cat\'egorie $\mathcal{P}(B)$ (avec les isomorphismes comme morphismes). D'apr\`es 1.7,
nous en d\'eduisons que les cat\'egories ${}_{\varepsilon}\mathcal{Q}^{\hbox{\tiny{\'el}}}(A)$ et 
$\mathcal{P}(B(e))$ sont \'equivalentes. En effet, nous avons
montr\'e en 1.7 que le groupe des automorphismes d'un objet de 
${}_{\varepsilon}\mathcal{Q}^{\hbox{\tiny{\'el}}} (A)$ est le m\^eme que celui
des automorphismes du $B$-module correspondant, vu comme un objet de $B(e)$ par extension des
scalaires. Puisque les classes d'isomorphie d'objets de $\mathcal{P}(B(e))$ co\"{i}ncident avec les classes
d'isomorphie d'objets de $\mathcal{P}(B)$, l'assertion r\'esulte de consid\'erations g\'en\'erales sur les
\'equivalences de cat\'egories.

\subsection{} Rappelons maintenant la d\'efinition du foncteur hyperbolique classique
$$H : \mathcal{P}(A) \to {}_{\varepsilon}\mathcal{Q}^{\hbox{\tiny{min}}}(A)$$
Si $E$ est un objet de $\mathcal{P}(A)$, $H(E)$ est le $A$-module $E \oplus E^*$ muni de la forme quadratique
$${\varphi}_0 : E \oplus E^* \to (E \oplus E^*)^* \approx E^*\oplus E$$
d\'efinie par la matrice
$$\varphi_0 =\left(
\begin{array}{cc}
0& 1\\
0& 0
\end{array}\right)$$
Si $u$ est un isomorphisme dans la cat\'egorie $\mathcal{P}(A)$, on d\'efinit 
$H(u) = g = u \oplus {}^tu^{-1}$. On v\'erifie
que ${}^tg.{\varphi}_0 .g = {\varphi}_0$ et que $H(u)$ est donc bien un isomorphisme dans la cat\'egorie 
${}_{\varepsilon}\mathcal{Q}^{\hbox{\tiny{min}}}(A)$.
On peut d\'ecrire ce foncteur de mani\`ere plus conceptuelle en consid\'erant l'anneau $\Lambda = M_2(A)$
des matrices $2\times 2$ \`a coefficients dans $A$ et o\`{u} l'involution est d\'efinie par
$$\left(
\begin{array}{cc}
a& b\\
c& d
\end{array}\right)
\mapsto 
\left(
\begin{array}{cc}
\overline{d}& \overline{b}\\
\overline{c}& \overline{a}
\end{array}\right)$$
L'\'equivalence de Morita d\'emontr\'ee dans \cite{B} \S  9 montre que les cat\'egories 
${}_{\varepsilon}\mathcal{Q}^{\hbox{\tiny{min}}}(\Lambda)$ et
${}_{\varepsilon}\mathcal{Q}^{\hbox{\tiny{min}}}(A)$ sont \'equivalentes. On d\'emontre par la m\^eme 
m\'ethode que les cat\'egories ${}_{\varepsilon}\mathcal{Q}^{\hbox{\tiny{max}}}(\Lambda)$ et
${}_{\varepsilon}\mathcal{Q}^{\hbox{\tiny{max}}}(A)$ d'une part et les cat\'egories 
${}_{\varepsilon}\mathcal{Q}^{\hbox{\tiny{\'el}}} (\Lambda)$ et 
${}_{\varepsilon}\mathcal{Q}^{\hbox{\tiny{\'el}}} (A)$ d'autre part sont \'equivalentes.
Le foncteur hyperbolique $\mathcal{P}(A)\to {}_{\varepsilon}\mathcal{Q}^{\hbox{\tiny{min}}}(A)$ 
est alors induit par l'homomorphisme d'anneaux $A\times A^{op} \to M_2(A)$ d\'efini par
$$(a, b)\mapsto \left(
\begin{array}{cc}
a& 0\\
0& \overline{b}
\end{array}\right)
$$
D'apr\`es 1.8, cette m\'ethode a l'avantage de d\'efinir un nouveau foncteur hyperbolique de $\mathcal{P}(A)$
dans la cat\'egorie plus fine ${}_{\varepsilon}\mathcal{Q}^{\hbox{\tiny{\'el}}}(A)$ par la composition 
des foncteurs \'evidents suivants induits par des morphismes d'anneaux ou des \'equivalences de Morita :
$$\mathcal{P}(A) \to \mathcal{P}(A(e)) {\sim} 
{}_{\varepsilon}\mathcal{Q}^{\hbox{\tiny{\'el}}} (A\times A^{op})\to
{}_{\varepsilon}\mathcal{Q}^{\hbox{\tiny{\'el}}} (M_2 (A)) {\sim} 
{}_{\varepsilon}\mathcal{Q}^{\hbox{\tiny{\'el}}} (A)$$
On proc\`ede de m\^eme pour le foncteur ``oubli'' 
${}_{\varepsilon}\mathcal{Q}^{\hbox{\tiny{\'el}}}(A)\to\mathcal{P}(A)$ qui est la composition
$${}_{\varepsilon}\mathcal{Q}^{\hbox{\tiny{\'el}}}(A) \to 
{}_{\varepsilon}\mathcal{Q}^{\hbox{\tiny{\'el}}} (A \times A^{op} ) {\sim}  \mathcal{P}(A(e)) \to \mathcal{P}(A)$$

\section{Les groupes de Grothendieck et Bass en $K-$th\'eorie hermitienne.}

\subsection{} Aux cat\'egories pr\'ec\'edentes 
${}_{\varepsilon}\mathcal{Q}^{\hbox{\tiny{max}}}(A)$, ${}_{\varepsilon}\mathcal{Q}^{\hbox{\tiny{min}}}(A)$ et 
${}_{\varepsilon}\mathcal{Q}^{\hbox{\tiny{\'el}}}(A)$, nous pouvons associer trois
groupes de $K-$th\'eorie hermitienne not\'es 
${}_{\varepsilon}K\mathcal{Q}^{\hbox{\tiny{max}}}(A)$, ${}_{\varepsilon}K\mathcal{Q}^{\hbox{\tiny{min}}}(A)$ 
et ${}_{\varepsilon}K\mathcal{Q}^{\hbox{\tiny{\'el}}}(A)$ respectivement, reli\'es par des homomorphismes canoniques
$$\xymatrix{{}_\varepsilon K\mathcal{Q}^{\hbox{{\small{\'el}}}}(A) \ar[r]^u&{}_\varepsilon K\mathcal{Q}^{
\operatorname{min}}(A)\ar[r]^v&{}_\varepsilon K\mathcal{Q}^{\operatorname{max}}(A)},$$
Il est clair que $u$ est un isomorphisme et que $v$ est surjectif. Par ailleurs, nous pouvons d\'efinir
l'analogue du groupes de Bass $K_1(A)$ en $K-$th\'eorie hermitienne. Dans ce but, le lemme suivant,
dont la d\'emonstration est d\'etaill\'ee dans \cite{KV} p. 61 par exemple, est essentiel.

{\lemma{Tout module ${\varepsilon}$-quadratique est facteur direct d'un module hyperbolique.}}

\subsection{} Puisque tout module projectif de type fini est facteur direct d'un module libre du type $A^n$,
on voit que les groupes classiques qui jouent le r\^ole de $GL_n(A)$ sont les groupes
d'automorphismes de modules hyperboliques du type $H(A^n)$ dans chacune des trois cat\'egories
concern\'ees.\\
Plus pr\'ecis\'ement, \'ecrivons $E = M \oplus M^*$ (on consid\`erera le cas o\`{u} $M=A^n$ un peu plus tard).
La forme quadratique associ\'ee est d\'efinie par la matrice ${\phi}_0$ pr\'ec\'edente avec $\psi$ comme forme
hermitienne associ\'ee, soit
$$\phi_0 = 
\left(\begin{array}{cc}
0& 1\\
0 &0
\end{array}\right)
\hbox{ }\phi=
\left(\begin{array}{cc}
0& 1\\
\varepsilon & 0
\end{array}\right)$$
Si $f:E\to E$ est un homomorphisme d\'efini par une matrice
$f =\left(\begin{array}{cc} 
a& b\\
c& d
\end{array}\right)$, 
son adjoint est la matrice $f^*=\left(\begin{array}{cc}
^td &\varepsilon ^tb\\
\varepsilon ^tc &^ta
\end{array}\right)$.\\
Dans le cas o\`{u} $M = A^n$, il convient de remplacer la notation ${}^tu$ par ${}^t\overline{u}$ , si on \'ecrit 
$u$ comme une matrice $n\times n$. En effet, la conjugaison r\'esulte de l'identification de $A^n$ avec 
son dual $(A^n)^*$.

\subsection{Notations} On d\'esigne par ${}_{\varepsilon}O_{n,n}^{\hbox{\tiny{max}}}(A)$ (resp. 
${}_{\varepsilon}O_{n,n}^{\hbox{\tiny{min}}}(A)$ , ${}_{\varepsilon}O_{n,n}^{\hbox{\tiny{\'el}}}(A)$) 
le groupe unitaire (resp. orthogonal, orthogonal \'elargi) associ\'e au module hyperbolique 
$(A)^n \oplus (A^n)^*$.

\subsection{Exemple} Supposons que $A$ soit un corps muni de l'involution triviale et que ${\varepsilon}=1$. 
Le fait que $f$ soit unitaire ($f \in  {}_1O_{n,n}^{\hbox{\tiny{max}}}(A)$) se traduit par les 
identit\'es suivantes (o\`{u} $a$, $b$, $c$ et $d$ sont des matrices $n\times n$) :
\begin{eqnarray*}
a.{}^td + b.{}^tc = 1\\
a.{}^tb + b.{}^ta = 0\\
c.{}^td + d.{}^tc = 0\\
c.{}^tb + d.{}^ta = 1
\end{eqnarray*}
L'automorphisme $f$ est orthogonal ($f \in  {}_1O_{n,n}^{\hbox{\tiny{min}}}(A)$) s'il existe en outre des 
matrices $h$ et $k$ telles
que $a.{}^tb = h - {}^th$ et $c.{}^td = k - {}^tk$.\\ \\
Pour d\'ecrire un \'el\'ement du groupe orthogonal \'elargi, il faut se donner en outre un
endomorphisme d\'efini par une matrice $2n\times 2n$
$$u =
\left(\begin{array}{cc}
\alpha & \beta\\
\gamma  &\delta
\end{array}\right)$$
li\'ee \`a $f$ et \`a la forme ${\phi}_0$ (cf. 1.6/7). Plus pr\'ecis\'ement, le couple $(f, u)$ 
doit v\'erifier l'identit\'e suivante
$$\left(\begin{array}{cc}
^td.a& ^tb.d\\
^tc.a& ^tc.b
\end{array}\right)=
\left(\begin{array}{cc}
1 &0\\
0 &0
\end{array}\right) +
\left(\begin{array}{cc}
\alpha & \beta \\
\gamma  &\delta 
\end{array}\right) -
\left(\begin{array}{cc}
^t\delta & ^t\beta \\
^t\gamma &^t \alpha 
\end{array}\right)$$
Elle r\'esulte de l'\'equation $(E)$ en 1.7, \`a condition d'identifier $E\oplus E^*$ \`a $E^*\oplus E$ 
(avec $E = A^n$).

\subsection{} Revenons au cas g\'en\'eral d'un anneau $A$ quelconque. Pour simplifier, nous \'ecrirons
${}_{\varepsilon}O_{n,n}(A)$ au lieu de ${}_{\varepsilon}O_{n,n}
^{\hbox{\tiny{max}}}(A)$, ${}_{\varepsilon}O_{n,n}
^{\hbox{\tiny{min}}}(A)$, ${}_{\varepsilon}O_{n,n}
^{\hbox{\tiny{\'el}}}(A)$ en revenant \`a ces notations sp\'ecifiques
lorsqu'il sera n\'ecessaire de distinguer les trois groupes. De m\^eme, nous utiliserons la
terminologie uniforme ``groupe orthogonal'' au lieu de ``groupe unitaire'', ``groupe orthogonal''
ou ``groupe orthogonal \'elargi'', lorsque nos consid\'erations s'appliquent aux trois variantes. Avec
ces conventions, le groupe orthogonal infini ${}_{\varepsilon}O(A)$ est d\'efini comme la limite inductive des
groupes ${}_{\varepsilon}O_{n,n}(A)$ avec les inclusions \'evidentes. En suivant l'exemple du groupe lin\'eaire, nous
d\'efinissons le ``groupe de Bass'' ${}_{\varepsilon}K\mathcal{Q}_{1}(A)$ comme le quotient de 
${}_{\varepsilon}O(A)$ par le sous-groupe des
commutateurs $[{}_{\varepsilon}O(A), {}_{\varepsilon}O(A)]$. Le fait que ce sous-groupe soit parfait 
r\'esulte de consid\'erations bien connues sur la stabilisation des matrices qu'on peut r\'esumer par 
des identit\'es g\'en\'erales. La premi\`ere est la suivante :
$$\tiny{\left(\begin{array}{ccc}
{\alpha}{\beta}{\alpha}^{-1}{\beta}& 0& 0\\
0& 1& 0\\
0& 0& 1
\end{array}\right)
=\left(
\begin{array}{ccc}
{\alpha}& 0& 0\\
0& {\alpha}^{-1}& 0\\
0 &0 &1
\end{array}\right)
\left(
\begin{array}{ccc}
{\beta}& 0& 0\\
0& 1& 0\\
0 &0 &{\beta}^{-1}
\end{array}\right)
\left(
\begin{array}{ccc}
{\alpha}^{-1}& 0& 0\\
0 &{\alpha}& 0\\
0 &0& 1
\end{array}\right)
\left(
\begin{array}{ccc}
{\beta}^{-1}& 0& 0\\
0& 1 &0\\
0 &0 &{\beta}
\end{array}\right)}$$
Par ailleurs, modulo le sous-groupe des commutateurs, une matrice du type
$$\left(
\begin{array}{ccc}
{\alpha}& 0& 0\\
0& {\alpha}^{-1}& 0\\
0 &0& 1
\end{array}\right)$$
peut aussi s'\'ecrire
$$\left(
\begin{array}{ccc}
{\alpha}& 0& 0\\
0& {\alpha}^{-1}& 0\\
0& 0& 1
\end{array}\right)
\left(
\begin{array}{ccc}
0& 1& 0\\
0& 0 &1\\
1& 0& 0
\end{array}\right)
=
\left(
\begin{array}{ccc}
0 &{\alpha}& 0\\
0 &0 &{\alpha}^{-1}\\
1 &0& 0
\end{array}\right)$$
qui est le commutateur suivant
$$\left(
\begin{array}{ccc}
{\alpha}& 0& 0\\
0 &0 &1\\
0 &1& 0
\end{array}\right)
\left(
\begin{array}{ccc}
0 &1& 0\\
1 &0& 0\\
0 &0& 1
\end{array}\right)
\left(
\begin{array}{ccc}
{\alpha}^{-1}& 0& 0\\
0& 0 &1\\
0 &1 &0
\end{array}\right)
\left(
\begin{array}{ccc}
0& 1& 0\\
1& 0& 0\\
0& 0& 1
\end{array}\right)$$
Toutes ces identit\'es (qui sont vraies dans le cadre plus g\'en\'eral de cat\'egorie mono\"{i}dales
sym\'etriques) d\'emontrent bien que $[{}_{\varepsilon}O(A), {}_{\varepsilon}O(A)]$ est parfait.
Pour chacune des trois th\'eories consid\'er\'ees, on utilisera les notations ${}_{\varepsilon} K\mathcal{Q}_1
^{\hbox{\tiny{max}}}(A)$, ${}_{\varepsilon} K\mathcal{Q}_1^{\hbox{\tiny{min}}}(A)$, 
${}_{\varepsilon} K\mathcal{Q}_1^{\hbox{\tiny{\'el}}}(A)$ ou simplement ${}_{\varepsilon}K\mathcal{Q}_{1}(A)$.

{\theorem{} {Consid\'erons un carr\'e cart\'esien d'anneaux hermitiens (avec $\varphi_1$ surjectif)
$$\xymatrix{A\ar[d]_{\varphi_2}\ar[r]^{\psi_1} & A_1\ar[d]^{\varphi_1}\\
A_2\ar[r]_{\psi_2}&A'}$$
On a alors une suite exacte (dite de Mayer-Vietoris) entre les groupes de $K$-th\'eorie
hermitienne
$$\begin{array}{ll}
\xymatrix{
{}_{\varepsilon}K\mathcal{Q}_1(A)\ar[r]& 
{}_{\varepsilon}K\mathcal{Q}_1(A_1)\oplus {}_{\varepsilon}K\mathcal{Q}_1(A_2)\ar[r]& 
{}_{\varepsilon}K\mathcal{Q}_1(A')\ar[r]& {}_{\varepsilon}K\mathcal{Q}(A)\ar[r]&\dots}\\
\qquad\xymatrix{\dots\ar[r]&  
{}_{\varepsilon}K\mathcal{Q}(A_1)\oplus {}_{\varepsilon}K\mathcal{Q}(A_2)\ar[r]&  {}_{\varepsilon}K\mathcal{Q}(A')}
\end{array}$$}}
\demo Ce th\'eor\`eme classique peut \^etre d\'emontr\'e de diverses mani\`eres. L'une d'entre
elle est esquiss\'ee dans le livre de Milnor \cite{M} et d\'etaill\'ee dans celui de Bak \cite{B}. Une autre
d\'emonstration est indiqu\'ee dans \cite{KV} p. 68-70 (elle s'applique dans les trois situations). Le
point important est de remarquer qu'un \'el\'ement du sous-groupe des commutateurs
$[{}_{\varepsilon}O(A'), {}_{\varepsilon}O(A')]$ se rel\`eve en un \'el\'ement de ${}_{\varepsilon}O(A_1)$. 
Ceci est d\'emontr\'e gr\^{a}ce au lemme de Whitehead classique adapt\'e au cas hermitien (cf. \cite{KV} 
th\'eor\`eme 2.6 par exemple).

\subsection{} Dans \cite{B} p. 191, Bak d\'emontre une suite exacte int\'eressante reliant les groupes 
${}_{\varepsilon}K\mathcal{Q}^{\hbox{\tiny{max}}}$ et ${}_{\varepsilon}K\mathcal{Q}^{\hbox{\tiny{min}}}$. 
Elle s'\'ecrit 
$$\xymatrix{{}_{\varepsilon} K\mathcal{Q}_1^{\hbox{\tiny{min}}}(A) \ar[r]& {}_{\varepsilon} K\mathcal{Q}_1
^{\hbox{\tiny{max}}}(A) \ar[r]&{}_{\varepsilon}{\Xi}(A) \ar[r]& {}_{\varepsilon} K\mathcal{Q}^{\hbox{\tiny{min}}}(A)
\ar[r]&{}_{\varepsilon} K\mathcal{Q}^{\hbox{\tiny{max}}}(A)}$$
Le groupe de 2-torsion ${}_{\varepsilon}{\Xi}(A)$ est explicit\'e ainsi. Nous d\'efinissons d'abord ${\Gamma}=\Gamma(A)$ comme
l'ensemble des \'el\'ements $a$ de $A$ tels que $\overline{a}={\varepsilon}a$ et $\Lambda$ comme le sous-groupe de 
${\Gamma}$ form\'e des $b-{\varepsilon}\overline{b}$. Alors ${}_{\varepsilon}{\Xi}(A)$ est le quotient de 
${\Gamma}/\Lambda {\otimes}_A{\Gamma}/\Lambda$ par le sous-groupe engendr\'e par tous les
\'el\'ements de la forme
$$\{a {\otimes} b - b {\otimes} a\}\hbox{ et }\{a {\otimes} b - a {\otimes} b a \overline{b} \}$$
Dans la d\'efinition du produit tensoriel ${\Gamma}/\Lambda{\otimes}_A{\Gamma}/\Lambda$, l'action \`a 
droite de $A$ sur ${\Gamma}/\Lambda$ est $({\gamma}, a){\mapsto}\overline{a} {\gamma} a$. L'action 
\`a gauche est d\'efinie de mani\`ere similaire par $(a, {\gamma}) {\mapsto} a .{\gamma}.\overline{a}$
Un th\'eor\`eme plus g\'en\'eral est en fait \'enonc\'e dans \cite{B} en utilisant des ``formes param\`etres''
arbitraires ${\Gamma}$ et $\Lambda$.

\subsection{Remarque.} La suite exacte pr\'ec\'edente permet de d\'efinir un invariant des formes
quadratiques proche de l'invariant de Arf en consid\'erant des corps de caract\'eristique 2 (cf.
[B2]). Dans ce cas, le groupe $K\mathcal{Q}_1^{\hbox{\tiny{max}}}(A)$ est r\'eduit \`a 0, $K\mathcal{Q}^{\hbox{\tiny{max}}} (A)\cong\mathbb{Z}$ et le noyau de la fl\`eche
$$K\mathcal{Q}^{\hbox{\tiny{min}}} (A) \to K\mathcal{Q}^{\hbox{\tiny{max}}} (A)=\mathbb{Z}$$
s'identifie ainsi au groupe ${\Xi}(A)$ pr\'ec\'edent : c'est le quotient de 
$A {\otimes}_{\mathbb{Z}}A$ par le sous-groupe
engendr\'e par les relations $\{a {\otimes} b - b {\otimes} a\}$, $\{a {\otimes} b - a {\otimes} b^2 a\}$ et 
$\{c^2a {\otimes} b - a {\otimes} c^2b\}$.
L'invariant de Arf classique est obtenu par l'application $a {\otimes} b {\mapsto} a.b$ : elle est \`a valeurs dans
le quotient $G$ de $F$ par le sous-groupe additif engendr\'e par les relations $\{a^2 - a\}$. Cette
application de ${\Xi}(A)$ dans G admet une r\'etraction induite par l'application $a {\mapsto} 1 {\otimes} a$.

\section{Les groupes de $K-$th\'eorie hermitienne ${}_{\varepsilon}K\mathcal{Q}_{n}(A)$ pour $n < 0$ et $n > 0$.}

\subsection{} Pour d\'efinir les groupes ${}_\varepsilon K\mathcal{Q}_{n}$ pour $n < 0$, nous suivons le m\^eme 
sch\'ema qu'en $K$-th\'eorie
alg\'ebrique \cite{KV}. De mani\`ere pr\'ecise, si on pose $n = -m$, on pose $ K\mathcal{Q}_{n}(A) = KQ(S^mA)$, o\`u
$S^mA$ est la $m^{\hbox{\tiny{i\`eme}}}$ suspension de l'anneau $A$. Notons que l'isomorphisme
$${}_\varepsilon K\mathcal{Q}^{\hbox{{\small{\'el}}}}(A) \cong   {}_\varepsilon K\mathcal{Q}^{\operatorname{min}}(A)$$ 
implique par suspensions it\'er\'ees l'isomorphisme
$${}_\varepsilon  K\mathcal{Q}_{-m}
^{\hbox{{\small{\'el}}}} (A) \cong   {}_\varepsilon  K\mathcal{Q}_{-m}^{\operatorname{min}}(A)$$ 
Le th\'eor\`eme suivant est moins \'evident.

{\theorem{}{L'homomorphisme
$${}_{\varepsilon} K\mathcal{Q}_n
^{\hbox{\tiny{min}}}(A)\to {}_{\varepsilon} K\mathcal{Q}_n
^{\hbox{\tiny{max}}}(A)$$
est surjectif pour $n=0$, bijectif pour $n<0$.}}
\demo La surjectivit\'e pour tout $n$ est une cons\'equence imm\'ediate des d\'efinitions (car
nous consid\'erons des formes hermitiennes paires). Par induction sur $n$, il suffit de d\'emontrer
l'injectivit\'e pour $n=-1$. Pour cela, \'ecrivons la suite exacte 2.8, en rempla\c{c}ant $A$ par sa
suspension $SA$ et son c\^one $CA$. On obtient alors un diagramme commutatif

$$\footnotesize{\xymatrix{{}_\varepsilon K\mathcal{Q}_1
^{\operatorname{min}}(CA) \ar[r]\ar[d]& {}_\varepsilon K\mathcal{Q}_1
^{\operatorname{max}}(CA) \ar[r]\ar[d]& {}_\varepsilon  \Xi (CA) \ar[r]\ar[d]& {}_\varepsilon K\mathcal{Q}^{
\operatorname{min}} (CA) \ar[r]\ar[d]& {}_\varepsilon K\mathcal{Q}^{\operatorname{max}} (CA)\ar[d]\\
{}_\varepsilon K\mathcal{Q}_1
^{\operatorname{min}}(SA)
\ar[r]&
{}_\varepsilon K\mathcal{Q}_1
^{\operatorname{max}}(SA) \ar[r]& {}_\varepsilon  \Xi (SA) \ar[r]& {}_\varepsilon K\mathcal{Q}^{\operatorname{min}} (
SA) \ar[r]&  {}_\varepsilon K\mathcal{Q}^{\operatorname{max}} (SA)}}$$
Puisque le c\^one d'un anneau est ``flasque'' (il existe un foncteur ${\tau}$ de la cat\'egorie 
$\mathcal{P}(CA)$ dans elle-m\^eme tel que ${\tau}  \oplus Id$ soit isomorphe \`a ${\tau}$ ), ses groupes de 
$K-$th\'eorie hermitienne sont r\'eduits \`a 0, ce qui implique que ${}_{\varepsilon}{\Xi}(CA)$ est aussi 
\'egal \`a 0. Pour d\'emontrer l'injectivit\'e de la fl\`eche ${}_{\varepsilon} K\mathcal{Q}_{-1}
^{\hbox{\tiny{min}}}(A)\to {}_{\varepsilon} K\mathcal{Q}_{-1}
^{\hbox{\tiny{max}}}(A)$, il suffit donc de montrer que l'homomorphisme ${}_{\varepsilon}{\Xi}(CA)
\to {}_{\varepsilon}{\Xi}(SA)$ est surjectif, ce qui est une cons\'equence du lemme suivant.

{\lemma{ Notons ${\Gamma}(R)$ le groupe ${\Gamma}$ d\'efini en 2.8 pour tout anneau $R$. 
Alors l'homomorphisme canonique
$${\Gamma}(CA) \to {\Gamma}(SA)$$
est surjectif.}}
\demo Un \'el\'ement de ${\Gamma}(SA)$ est d\'efini par une matrice infinie $M$ telle que sur chaque
ligne et chaque colonne il n'existe qu'un nombre fini d'\'el\'ements non nuls et telle que ${}^t\overline{M} = {\varepsilon} M$
modulo une matrice finie. Soient $a_{ij}$ les \'el\'ements (en nombre fini) de la matrice $M$ tels que 
$\overline{a_{ij}} \neq {\varepsilon} a_{ji}$ . Si on remplace ces \'el\'ements par 0, on trouve une matrice $N$ dans $CA$ 
qui est ${\varepsilon}$-hermitienne et dont la classe dans $SA$ est \'egale \`a celle de $M$.

\subsection{} D\'efinissons maintenant les groupes ${}_{\varepsilon}K\mathcal{Q}_{n}$ pour $n > 0$, ce qui est 
plus d\'elicat. En principe, il suffit de copier la construction + de Quillen \`a l'espace $B{}_{\varepsilon}O(A)$, 
ce qui est possible car le sous-groupe des commutateurs $[{}_{\varepsilon}O(A), {}_{\varepsilon}O(A)]$ est parfait. 
On d\'efinit alors ${}_{\varepsilon}K\mathcal{Q}_{n}(A)$ comme le $n^{\hbox{\tiny{i\`eme}}}$ 
groupe d'homotopie de $B{}_{\varepsilon}O(A)^+$ (pour $n > 0$). En fait, nous disposons de trois groupes de 
$K$-th\'eorie hermitienne
$${}_{\varepsilon} K\mathcal{Q}_n
^{\hbox{\tiny{max}}}(A),\hbox{ }{}_{\varepsilon} K\mathcal{Q}_n
^{\hbox{\tiny{min}}}(A)\hbox{ et }{}_{\varepsilon} K\mathcal{Q}_n
^{\hbox{\tiny{\'el}}}(A)$$
associ\'es respectivement aux groupes ${}_{\varepsilon}O^{\hbox{\tiny{max}}}(A)$, 
${}_{\varepsilon}O^{\hbox{\tiny{min}}}(A)$ et ${}_{\varepsilon}O^{\hbox{\tiny{\'el}}}(A)$. Conform\'ement \`a la
philosophie de cet article, nous adopterons la notation uniforme ${}_{\varepsilon}K\mathcal{Q}_{n}(A)$ pour ne pas
compliquer l'exposition, lorsqu'il n'y a pas de risque de confusion. Ces groupes sont difficiles \`a
calculer en g\'en\'eral, comme d'ailleurs les groupes $K_n(A)$ de Quillen dont ils sont la
g\'en\'eralisation. Nous verrons cependant que, dans une certaine mesure, les ``groupes de Witt
sup\'erieurs'' ${}_{\varepsilon}W_n(A) = Coker(K_n(A) \to {}_{\varepsilon}K\mathcal{Q}_{n}(A))$ sont plus accessibles.

\subsection{} Comme il est bien connu, il existe d'autres d\'efinitions des foncteurs $K_n$ et 
${}_{\varepsilon}K\mathcal{Q}_{n}$ \'equivalentes \`a la construction + de Quillen. La construction dite 
``$S^{-1}S$'' (due aussi \`a Quillen) est d\'etaill\'ee dans le cadre hermitien dans \cite{K1} \S  1 et nous 
l'utiliserons pour la preuve de 4.5. Il existe aussi une d\'efinition en termes de $A$-fibr\'es plats qui 
est d\'etaill\'ee dans \cite{K2} p. 42 et c'est celle que nous utiliserons essentiellement ici. Rappelons-l\`a 
bri\`evement dans le cadre que nous int\'eresse.\\
On d\'efinit un $A$-fibr\'e hermitien ``virtuel'' sur un $CW$-complexe $X$ comme la donn\'ee d'une
fibration acyclique $Y \to X$ et d'un $A$-fibr\'e plat $E$ sur $Y$, la fibre \'etant un $A$-module projectif de
type fini muni d'une forme hermitienne dans l'un des trois sens que nous avons donn\'es \`a ce
terme (ceci veut dire que les fonctions de transition du fibr\'e sur $Y$ sont des fonctions localement
constantes dans chacune des trois cat\'egories ``max'', ``min'' ou ``\'el'' concern\'ees).\\ \\
Deux tels fibr\'es virtuels
$$E \to Y \to X\hbox{ et }E' \to Y' \to X$$
sont dits \'equivalents s'il existe un fibr\'e virtuel $E_1 \to Y_1 \to X$ et un diagramme commutatif
$$\xymatrix{Y \ar[r]^f\ar[d]_\sigma& X\ar[d]^{f'}\\
Y_1\ar[ru]^{f_1}&Y'\ar[l]^{\sigma'}}$$
tel que $\sigma^*(E_1 ) \cong   E$ et $\sigma'^*(E_1 ) \cong   E'$.\\

En suivant le m\^eme sch\'ema qu'en \cite{K2} p. 42-50, on montre que le groupe de Grothendieck
construit avec ces fibr\'es virtuels est isomorphe au groupe d\'efini par les classes d'homotopie de
$X$ dans ${}_{\varepsilon}K\mathcal{Q}_0(A)\times B{}_{\varepsilon}O(A)^+$, not\'e 
${}_{\varepsilon}K\mathcal{Q}_A(X)$, et qui est une ``th\'eorie cohomologique'' en $X$. Si $X$
est une sph\`ere de dimension $n \geq  0$, on retrouve ainsi ${}_{\varepsilon}K\mathcal{Q}_{n}(A)$ 
comme le conoyau de la fl\`eche \'evidente ${}_{\varepsilon}K\mathcal{Q}_0(A)\to{}_{\varepsilon}K\mathcal{Q}_A(X)$.
\\ \\
On peut d\'efinir le spectre de la $K-$th\'eorie hermitienne par la m\^eme m\'ethode qu'en $K-$th\'eorie
alg\'ebrique. Ainsi, dans \cite{K1}, on d\'emontre l'analogue du th\'eor\`eme de Gersten-Wagoner \cite{W} en
$K-$th\'eorie hermitienne : on a une \'equivalence d'homotopie (non naturelle) entre ${\Omega}(B{}_{\varepsilon}O(SA)^+)$ et
${}_{\varepsilon}K\mathcal{Q}_0(A)\times B{}_{\varepsilon}O(A)^+$ (la m\^eme d\'emonstration
 s'applique dans les trois cas consid\'er\'es ici). 
Plus pr\'ecis\'ement, on d\'efinit le $\Omega$-spectre 
de la $K$-th\'eorie hermitienne
$\bf{{}_{\varepsilon}KQ(A)}_*$ par les formules suivantes :
$$\begin{array}{lc}
\bf{{}_{\varepsilon}KQ}(A)_n = \Omega (B{}_{\varepsilon}O(S^{n+1}A)^+) &\hbox{pour $n\geq 0$}\\
\bf{{}_{\varepsilon}KQ}(A)_{n}= \Omega^{-n}(B{}_{\varepsilon}O(A)^+) &\hbox{pour $n < 0$}\end{array}$$
En fait, ce spectre n'est qu'un langage commode. Pour pouvoir d\'efinir des 
cup-produits en $K$-th\'eorie hermitienne, nous nous servirons plut\^ot de la th\'eorie cohomologique 
associ\'ee en termes de fibr\'es virtuellement plats comme nous l'avons explicit\'e plus haut. D'ailleurs, une
situation analogue se pr\'esente en $K-$th\'eorie topologique, o\`{u} les op\'erations sont
plus ais\'ement d\'efinies sur les fibr\'es vectoriels plut\^ot que sur la grassmannienne infinie.

\section{Cup-produits en $K-$th\'eorie hermitienne. Le cup-produit de Clauwens.}

\subsection{} L'avantage du point de vue des fibr\'es plats est une d\'efinition tr\`es simple du cup-produit.
Celui-ci est explicit\'e dans \cite{K2} \`a partir d'un morphisme $\mathbb{Z}$-bilin\'eaire

$$\varphi  : A \times B \to  C$$
v\'erifiant la propri\'et\'e de multiplicativit\'e suivante
$$\varphi (aa', bb') = \varphi (a, b) \varphi (a', b')$$
Le cup-produit s'\'ecrit alors sous la forme d'un accouplement bilin\'eaire
$$K_A(X) \times K_B(Y) \to  K_{A\otimes B}(X \times Y)$$
o\`u la fl\`eche est simplement induite par le produit tensoriel des fibr\'es virtuellement plats. Si $X$
est un espace muni d'un point base $P$, il est commode d'introduire la 
``$K$-th\'eorie r\'eduite" 
$\widetilde{K}_A(X)=Ker[K_A(X) \to  K_A(P)=K_0(A)]$. Le produit pr\'ec\'edent induit 
alors un ``cup-produit r\'eduit"
$$\widetilde K_A (X) \times \widetilde K_B (Y ) \to 
\widetilde K_{A\otimes B} (X\wedge Y )$$
En particulier, si $X$ (resp. $Y$) est une sph\`ere $S^n$ (resp. $S^p$) avec 
$n$ et $p \geq  0$, on en d\'eduit le cup-produit usuel en 
$K$-th\'eorie alg\'ebrique (cf. aussi \cite{L}).

\subsection{} Le m\^eme sch\'ema s'applique en $K-$th\'eorie hermitienne\footnote{\'A condition de supposer en outre que $\overline{{\varphi}(a,b)}=\varphi(\overline{a},\overline{b})$}. Par exemple, compte tenu des
signes de sym\'etrie, les cup-produits classiques sont sch\'ematis\'es par des accouplements
$${}_\varepsilon K\mathcal{Q}^{\operatorname{max}} \times {}_\eta K\mathcal{Q}^{\operatorname{min}} \to{}_{
\varepsilon \eta}K\mathcal{Q}^{\operatorname{min}}$$
et
$${}_\varepsilon K\mathcal{Q}^{\operatorname{max}} \times {}_\eta K\mathcal{Q}^{\hbox{{\small{\'el}}}} \to 
{}_{\varepsilon \eta}  K\mathcal{Q}^{\hbox{{\small{\'el}}}}$$
De mani\`ere pr\'ecise, si nous consid\'erons une $\varepsilon$-forme 
hermitienne paire $\phi  = \phi_0 + \varepsilon {}^t\phi_0$ sur un $A$-module
$E$ et une forme $\eta$-quadratique d\'efinie par une classe de de morphismes $\psi_0$ 
sur un $B$-module $F$, alors $\phi\otimes\psi_0$ est une classe de forme 
$\varepsilon \eta$-quadratique sur $E \otimes   F$. En outre, si $\alpha$ (resp. $\beta$ ) est un morphisme
unitaire (resp. orthogonal) de $E$ (resp. $F$), il est facile de voir que $\alpha\otimes\beta$ 
est un morphisme orthogonal de $E \otimes F$. De mani\`ere analogue, si $(\beta , \gamma)$ 
est un morphisme dans la cat\'egorie ${}_\eta\mathcal{Q}^{\hbox{{\small{\'el}}}}$, le couple
$(\alpha\otimes\beta,\alpha\otimes\gamma)$ d\'efinit un 
morphisme dans la cat\'egorie 
${}_{\varepsilon \eta}\mathcal{Q}^{\hbox{{\small{\'el}}}}$, ce qui d\'efinit 
le deuxi\`eme accouplement.\\ Ces deux cup-produits, d\'efinis en termes de modules, 
s'\'etendent naturellement aux fibr\'es plats ou virtuellement plats dans les cat\'egories concern\'ees 
(il convient de noter cependant que ${\psi}_0$ n'est pas donn\'e dans la structure pour le premier 
accouplement mais seulement sa classe fibre par fibre). En consid\'erant des fibr\'es plats sur 
des sph\`eres homologiques, on d\'efinit ainsi des accouplements
$${}_{\varepsilon} K\mathcal{Q}_n
^{\hbox{\tiny{max}}}(A)\times {}_{\eta} K\mathcal{Q}_p
^{\hbox{\tiny{min}}}(B)\to {}_{\varepsilon\eta} K\mathcal{Q}_{n+ p}
^{\hbox{\tiny{min}}}(C)$$
 et 
$${}_{\varepsilon} K\mathcal{Q}_n
^{\hbox{\tiny{max}}}(A)\times {}_{\eta} K\mathcal{Q}_p
^{\hbox{\tiny{\'el}}}(B)\to {}_{\varepsilon\eta} K\mathcal{Q}_{n+p}
^{\hbox{\tiny{\'el}}}(C)$$

\subsection{} Nous allons maintenant introduire un autre cup-produit plus subtil, d\^{u} essentiellement \`a
Clauwens \cite{C}. Celui-ci a \'et\'e \'ecrit par Clauwens pour les cat\'egories de modules mais il s'\'etend
ais\'ement aux ``bonnes'' cat\'egories des fibr\'es virtuellement plats munis de formes quadratiques.
De mani\`ere pr\'ecise, consid\'erons la cat\'egorie ${}_{\eta}\mathcal{Q}^{\hbox{\tiny{\'el}}}(B)$ 
ainsi que la sous-cat\'egorie ${}_{\varepsilon} \mathcal{Q}^{\hbox{\tiny{\'el}}_0}(A[s])$ de 
${}_{\varepsilon}\mathcal{Q}^{\hbox{\tiny{\'el}}}(A[s])$ form\'ee des 
$A[s]$-modules provenant de $A$ par extension des scalaires, 
l'involution sur $A[s]$ \'etant induite par l'involution de $A$ et la transformation $s {\mapsto} 1 - s$.\\
Un objet de ${}_{\varepsilon} \mathcal{Q}^{\hbox{\tiny{\'el}}}(A[s])$ peut \^etre d\'ecrit comme un couple 
$(E, {\theta})$, o\`{u} $E$ est un objet de $\mathcal{P}(A)$ et ${\theta}$
une forme ${\varepsilon}$-quadratique sur $E {\otimes}_\mathbb{Z}\mathbb{Z}[s]$ s'\'ecrivant sous la forme 
$\sum  {\theta}_n s^n$, o\`{u} ${\theta}_n$ est un morphisme de $E$ vers $E^*$.\\ 
Consid\'erons maintenant un objet $(F, {\delta})$ de ${}_{\eta}\mathcal{Q}^{\hbox{\tiny{\'el}}}(B)$ , o\`{u} 
${\delta}$ est une forme ${\eta}$-quadratique non d\'eg\'en\'er\'ee sur $F$ avec 
${\Delta}  = {\delta} + {\eta}{}^t{\delta}$ comme forme hermitienne associ\'ee. Sur $E {\otimes} F$ on peut
alors consid\'erer la forme ${\varepsilon}{\eta}$-quadratique d\'efinie par la formule suivante
$${\kappa}  = \sum  {\theta}_n {\otimes} {\Delta}  ({\Delta} ^{-1}{\delta})^n$$
Cette formule se simplifie si on identifie $F$ et son dual par l'isomorphisme ${\Delta}$, ce qui revient \`a
remplacer ${\Delta} ^{-1}{\delta}$ par ${\delta}$. On peut de m\^eme identifier 
$E$ \`a $E^*$ par l'isomorphisme ${\theta}_0 +\sum_{n=0}^{\infty}{}^t{\theta}_n$. Le
foncteur de dualit\'e $f {\mapsto} {}^{t}f$ est alors remplac\'ee par le foncteur d'adjonction $f {\mapsto} f^*$.
Un avantage de cette formulation est aussi de se d\'ebarrasser des signes de sym\'etrie.
La formule pr\'ec\'edente s'\'ecrit alors sous une forme plus simple
$${\kappa}  = \sum  {\theta}_n {\otimes} {\delta}^n$$
avec ${\delta}^* = 1 - {\delta}$.
En quelques lemmes fondamentaux (cf. \cite{C} p. 43 et 44 et aussi l'appendice, o\`{u} on \'ecrit ${\phi}$  au
lieu de ${\Delta}^{-1}{\delta}$ pour \'eviter toute confusion), Clauwens montre que l'accouplement pr\'ec\'edent
$$Obj ({}_{\varepsilon} \mathcal{Q}^{\hbox{\tiny{\'el}}_0} (A[s]) ) \times 
Obj ({}_{{\eta}}\mathcal{Q}^{\hbox{\tiny{\'el}}} (B) ) \to Obj ({}_{{\varepsilon}{\eta}}\mathcal{Q}^{\hbox{\tiny{\'el}}}(A{\otimes} B) )$$
est bien d\'efini sur les classes d'isomorphie de modules quadratiques \'elargis. En fait, Clauwens
consid\`ere dans son article des modules libres mais sa m\'ethode est plus g\'en\'erale, comme nous
l'explicitons dans l'appendice. En particulier, nous pouvons d\'efinir un cup-produit remarquable
$${}_{\varepsilon}K\mathcal{Q}^{\hbox{\tiny{\'el}}_0}  (A[s])\times{}_{\eta} 
K\mathcal{Q}^{\hbox{\tiny{\'el}}} (B) \to 
{}_{{\varepsilon}{\eta}} K\mathcal{Q}^{\hbox{\tiny{\'el}}} (A{\otimes} B)$$
o\`{u} ${}_{\varepsilon}K\mathcal{Q}^{\hbox{\tiny{\'el}}_0}$  ($A[s]$) est le sous-groupe de 
${}_{\varepsilon}K\mathcal{Q}^{\hbox{\tiny{\'el}}} (A[s])$ engendr\'e par les modules provenant de $A$
par extension des scalaires (ceci est stablement le cas si $A$ est noeth\'erien r\'egulier par exemple).\\ \\
Dans les consid\'erations pr\'ec\'edentes, nous aurions pu remplacer la cat\'egorie 
$\mathcal{Q}^{\hbox{\tiny{\'el}}}$ par la cat\'egorie plus simple $\mathcal{Q}^{\hbox{\tiny{min}}}$. 
La raison pour travailler dans la cat\'egorie $\mathcal{Q}^{\hbox{\tiny{\'el}}}$ est notre souhait
de g\'en\'eraliser l'acccouplement d\'efini sur les groupes $K\mathcal{Q}_0$ aux 
groupes $K\mathcal{Q}_{n}$ d\'efinis dans le \S  3 pour $n > 0$. Si nous choisissons la d\'efinition de la 
$K-$th\'eorie hermitienne en termes de fibr\'es plats, il nous faut montrer par exemple que la classe 
d'isomorphie de la forme quadratique ${\kappa}$ 
d\'efinie plus haut ne d\'epend que des classes de ${\theta}$ et de ${\delta}$. Les lemmes de Clauwens 
(red\'emontr\'es en appendice) montrent la n\'ecessit\'e de se donner le morphisme ${\gamma}$ 
dans la formule $(S)$ en 1.4. Gr\^{a}ce \`a ce nouveau point du vue, on peut \'etendre le cup-produit 
pr\'ec\'edent aux groupes de $K\mathcal{Q}$-th\'eorie sup\'erieurs (dans la cat\'egorie ``\'el''), soit
$${}_{\varepsilon} K\mathcal{Q}_{n}^{\hbox{\tiny{\'el}}_0}(A[s]) \times {}_{\eta} K\mathcal{Q}_p
^{\hbox{\tiny{\'el}}}(B) \to {}_{{\varepsilon}{\eta}}K\mathcal{Q}_{n+p}
^{\hbox{\tiny{\'el}}}(A{\otimes} B)$$

\subsection{} Au d\'ebut de son article (th\'eor\`eme 1, p. 42), Clauwens montre que modulo l'addition de $A$-modules
hyperboliques (voir l'appendice pour un \'enonc\'e pr\'ecis), on peut se ramener au cas o\`{u}
${\theta}$ est ``lin\'eaire'', i.e. du type ${\theta} = g s$. En d'autres termes, ${\theta}_n = 0$, 
\`a l'exception de ${\theta}_1$ qui est \'egal \`a $g$. Puisque la forme hermitienne associ\'ee 
$g s + {\varepsilon}{}^tg (1 - s)$ est un isomorphisme, ceci implique
que ${}^tg = {\varepsilon} g (1 + N)$, o\`{u} $N$ est un endomorphisme nilpotent de $E$ (un tel $g$ est dit ``presque
hermitien''). Dans ce cas, la formule pour la forme quadratique ${\kappa}$ ci-dessus est tr\`es simple : on
trouve 
$${\kappa}  = g {\otimes} {\delta}$$
(si on identifie $F$ \`a son dual par ${\Delta}$ )
En d'autres termes, l'accouplement pr\'ec\'edent sur les groupes $K\mathcal{Q}^{\hbox{\tiny{\'el}}}$ 
g\'en\'eralise (pour $N = 0$) l'accouplement classique entre les formes hermitiennes 
(non n\'ecessairement paires) et les formes quadratiques. Un cas particulier important est le cup-produit
$${}_1K\mathcal{Q}_1^{\hbox{\tiny{\'el}}_0}(S\mathbb{Z}[s]) \times {}_{\eta} K\mathcal{Q}_p
^{\hbox{\tiny{\'el}}}(B) \to{}_{\eta} K\mathcal{Q}_{1+p}
^{\hbox{\tiny{\'el}}}(S\mathbb{Z} {\otimes} B) = {}_{\eta} K\mathcal{Q}_{1+p}
^{\hbox{\tiny{\'el}}}(SB)$$

{\theorem{}{ Soit $u_1$ l'\'element de ${}_1K\mathcal{Q}_1^{\hbox{\tiny{\'el}}_0}(S\mathbb{Z}[s]) 
= {}_1K\mathcal{Q}_1^{\hbox{\tiny{\'el}}}(S\mathbb{Z}[s])$ correspondant \`a
l'\'el\'ement unit\'e dans $K\mathcal{Q}_0(\mathbb{Z}[s]) = \mathbb{Z}\times\mathbb{Z}$ (cf. \cite{C}, p. 47). 
Alors le cup-produit par $u_1$ induit un isomorphisme entre ${}_{\eta} K\mathcal{Q}_p^{\hbox{\tiny{\'el}}}(B)$ et 
${}_{\eta} K\mathcal{Q}_{1+ p}^{\hbox{\tiny{\'el}}}(SB)$}}
\demo Elle est analogue \`a celle en $K-$th\'eorie alg\'ebrique ou hermitienne classique (cf.
\cite{K1} p. 224).

\subsection{} Rappelons par ailleurs qu'un autre cup-produit plus simple a \'et\'e d\'efini en 4.2:
$${}_{\varepsilon} K\mathcal{Q}_{n}
^{\hbox{\tiny{max}}}(A) \times{}_{\eta} K\mathcal{Q}_p
^{\hbox{\tiny{\'el}}}(B) \to {}_{{\varepsilon}{\eta}} K\mathcal{Q}_{n+p}
^{\hbox{\tiny{\'el}}}(A{\otimes} B)$$
Ces deux produits sont reli\'es ainsi:

{\theorem{}{ Le cup-produit de Clauwens est partiellement associatif dans le sens
suivant. Pour trois anneaux $B$, $C$ et $D$, on a le diagramme commutatif (avec 
$n = n_1 + n_2$, $\varepsilon=\varepsilon_1 \varepsilon_2$)
$$\xymatrix{
{}_{\varepsilon_1}K\mathcal{Q}_{n_1}
^{\operatorname{max}}(C)\times{}_{\varepsilon_2} K\mathcal{Q}_{n_2}
^{\hbox{{\small{\'el}}}_0} (D[s])\times{}_\eta K\mathcal{Q}_{p}
^{\hbox{{\small{\'el}}}}(B) \ar[r]\ar[d]&{}_{\varepsilon_1} K\mathcal{Q}_{n_1}
^{\operatorname{max}}(C)\times {}_{\varepsilon_2 \eta}  K\mathcal{Q}_{n_2+p}
^{\hbox{{\small{\'el}}}} (D \otimes   B)\ar[d]\\
{}_{\varepsilon_1\varepsilon_2} K\mathcal{Q}_{n_1 +n_2}
^{\hbox{{\small{\'el}}}_0} ((C \otimes   D)[s]) \times{}_\eta  K\mathcal{Q}_{p}
^{\hbox{{\small{\'el}}}}(B) \ar[r]&{}_{\varepsilon\eta} K\mathcal{Q}_{n+p}
^{\hbox{{\small{\'el}}}} (C \otimes D \otimes B)
}$$}}
\demo C'est une cons\'equence directe de la formule donn\'ee en 4.3. Nous devons
multiplier les deux membres de la formule par la m\^eme 
forme hermitienne paire avant et apr\`es
avoir fait le produit tensoriel par $\Delta  (\Delta{}^{-1}\delta )^n$.

\remark Pour les degr\'es n\'egatifs, nous avons seulement \`a consid\'erer des modules sur
des suspensions it\'er\'ees des anneaux consid\'er\'es. La notion de forme quadratique \'elargie est
alors inutile dans les d\'emonstrations. On peut m\^eme se limiter aux formes hermitiennes paires
pour les degr\'es $<0$ d'apr\`es 3.2.

\remark Si 1 est scind\'e dans $A$ (par exemple si 2 est inversible), on 
a des isomorphismes
$K\mathcal{Q}_{n}
^{\hbox{{\small{\'el}}}}(A) \cong   K\mathcal{Q}_{n}
^{\operatorname{max}}(A(e)) \cong   K\mathcal{Q}_{n}
^{\operatorname{min}}(A(e))$ avec $\overline{e} = -e$.

\section{Le th\'eor\`eme fondamental de la $K-$th\'eorie hermitienne pour des anneaux arbitraires}

\subsection{} Dans ce paragraphe, nous allons d\'esigner le 
spectre de la $K$-th\'eorie hermitienne ainsi que
celui de la $K$-th\'eorie alg\'ebrique par des caract\`eres gras.
De mani\`ere pr\'ecise, $\bf{K(A)}$ 
repr\'esentera le spectre de la $K$-th\'eorie alg\'ebrique usuelle ; celui de la $K$-th\'eorie hermitienne sera repr\'esent\'e 
par l'un des trois
spectres $\bf{{}_{\varepsilon}K\mathcal{Q}^{\operatorname{max}}(A)}$, 
$\bf{{}_{\varepsilon}K\mathcal{Q}^{\operatorname{min}}(A)}$ ou 
$\bf{{}_{\varepsilon}K\mathcal{Q}^{\hbox{{\small{\'el}}}}(A)}$, suivant la th\'eorie 
consid\'er\'ee. En particulier, les
foncteurs ``oubli" et hyperbolique induisent des morphismes
$$\bf{{}_\varepsilon K\mathcal{Q}^{\hbox{{\small{\'el}}}}(A) \to  K(A)}\hbox{ et }
\bf{K(A) \to  {}_{-\varepsilon} K\mathcal{Q}^{\hbox{{\small{\'el}}}}(A)}$$
dont les fibres homotopiques respectives seront not\'ees $\bf{{}_\varepsilon  V^{\hbox{{\small{\'el}}}}(A)}$ et $\bf
{{}_{-\varepsilon} U^{\hbox{{\small{\'el}}}}(A)}$.
L'\'enonc\'e suivant g\'en\'eralise le th\'eor\`eme de \cite{K2} (p. 260).
{\theorem{}{Nous avons une \'equivalence d'homotopie naturelle
$$\bf{{}_\varepsilon V^{\hbox{{\small{\'el}}}}(A) \approx  \Omega{}_{-\varepsilon} 
U^{\hbox{{\small{\'el}}}}(A)}$$}}

\subsection{Remarques.} Le th\'eor\`eme est \'evident lorsque $A=B\times B^{op}$, une situation d\'ej\`a consid\'er\'ee
dans les paragraphes pr\'ec\'edents. Dans ce cas, les spectres 
$\bf{{}_\varepsilon V^{\hbox{{\small{\'el}}}}(A)}$ et 
$\bf{\Omega{}_{-\varepsilon} U ^{\hbox{{\small{\'el}}}}(A)}$ co\"{\i}ncident
tous les deux avec la fibre homotopique du morphisme \'evident $\bf{K(B(e))\to K(B) x K(B)}$.\\ \\
Par ailleurs, si 1 est scind\'e dans $A$, nous retrouvons le th\'eor\`eme fondamental de la $K$-th\'eorie hermitienne 
\'enonc\'e dans \cite{K2} p. 260 (cf. la
remarque 5.11 un peu plus loin). La d\'emonstration du th\'eor\`eme 5.2 va 
\^etre en fait calqu\'ee sur celle de \cite{K2}. Nous mentionnerons simplement ici 
les modifications \`a y apporter.

\subsection{} Rappelons d'abord le principe g\'en\'eral de la d\'emonstration dans \cite{K2}
que nous appliquerons \`a plusieurs reprises : un morphisme d'anneaux hermitiens 
$f:A\to B$ induit une application entre spectres
$$\bf{{}_{\varepsilon}K\mathcal{Q}^{\hbox{{\small{\'el}}}}(A)\to  {}_{\varepsilon}
K\mathcal{Q}^{\hbox{{\small{\'el}}}}(B)}$$
dont nous pouvons interpr\'eter la fibre homotopique 
d'apr\`es un argument adapt\'e de Wagoner
\cite{W}. Pour cela, on consid\`ere le produit fibr\'e d'anneaux
$$\xymatrix{
R \ar[r]\ar[d]&CB\ar[d]\\
SA \ar[r]&  SB}$$
d'o\`u on d\'eduit la fibration homotopique
$$\xymatrix{
\bf{{}_{\varepsilon}K\mathcal{Q}^{\hbox{{\small{\'el}}}}(R)} \ar[r]& \bf{{}_{\varepsilon}K\mathcal{Q}^{\hbox{{\small{\'el}}}}SA)} \ar[r]& \bf{{}_{\varepsilon}K\mathcal{Q}^{\hbox{{\small{\'el}}}}(SB)}}$$
car $\bf{{}_{\varepsilon}K\mathcal{Q} ^{\hbox{{\small{\'el}}}}(CB)}$ est contractile. L'espace des lacets de $\bf{{}_{\varepsilon}K\mathcal{Q}^{\hbox{{\small{\'el}}}}(R)}$ 
est donc la fibre homotopique
recherch\'ee du morphisme 
$$\bf{{}_{\varepsilon}K\mathcal{Q}^{\hbox{{\small{\'el}}}}(A)\to  
{}_{\varepsilon}K\mathcal{Q}^{\hbox{{\small{\'el}}}}(B)}$$

Deux cas importants peuvent \^etre consid\'er\'es. Dans le premier, le morphisme est 
$A \times A^{op}\to  M_2(A)$
et dans le second $A \to  A\times A^{op}$, tous les deux d\'efinis en 1.8. Si nous d\'esignons\footnote{En fait, 
pour la $K-$th\'eorie, c'est \`a dire la $K-$th\'eorie hermitienne de $A \times A^{op}$, nous devons 
remplacer l'anneau des nombres duaux $A(e)$ par $A$, comme il a \'et\'e pr\'ecis\'e en 1.8.} 
par $U_A$ (resp. $V_A$) l'anneau $R$ obtenu dans ces deux cas, nous voyons que 
$\bf{{}_{\varepsilon}U^{\hbox{{\small{\'el}}}}(A)}$ est homotopiquement
\'equivalent \`a $\bf{\Omega {}_{\varepsilon}K\mathcal{Q}^{\hbox{{\small{\'el}}}}(U_A)}$ et que 
$\bf{{}_{\varepsilon}V^{\hbox{{\small{\'el}}}}(A)}$ est homotopiquement \'equivalent \`a 
$\bf{\Omega K{}_{\varepsilon}\mathcal{Q}^{\hbox{{\small{\'el}}}}(V_A)}$.

\subsection{} Nous souhaitons d\'efinir une application
$$\bf{{}_\varepsilon V^{\hbox{{\small{\'el}}}}(SA) \to  {}_{-\varepsilon} U^{\hbox{{\small{\'el}}}}(A)}$$
L'id\'ee, d\'ej\`a pr\'esente dans \cite{K2}, est d'inclure cette application dans le diagramme suivant
$$\footnotesize{\xymatrix{
\bf{{}_\varepsilon K\mathcal{Q}^{\hbox{{\small{\'el}}}}(A)}\ar[r]\ar[d]&
\bf{K(A)}\ar[r]\ar[d]&
\bf{{}_\varepsilon V^{\hbox{{\small{\'el}}}}(SA)}\ar[r]\ar[d]&
\bf{{}_\varepsilon K\mathcal{Q}^{\hbox{{\small{\'el}}}}(SA)}\ar[r]\ar[d]^\sigma&
\bf{K(SA)}\ar[r]\ar[d]&\bf{{}_\varepsilon V^{\hbox{{\small{\'el}}}}(S^2A)}\ar[d]\\
\bf{{}_{-\varepsilon} D^{\hbox{{\small{\'el}}}}(A)}\ar[r]&
\bf{K(A)}\ar[r]&\bf{{}_{-\varepsilon} U^{\hbox{{\small{\'el}}}}(A)}\ar[r]&\bf{{}_{-
\varepsilon}D^{\hbox{{\small{\'el}}}}(SA)}\ar[r]&\bf{K(SA)}
\ar[r]&\bf{{}_{-\varepsilon} U^{\hbox{{\small{\'el}}}}(SA)
}}}$$

La th\'eorie $\bf{{}_{-\varepsilon} D^{\hbox{{\small{\'el}}}}(A)}$ est ici la 
fibre homotopique de l'application $\hbox{\small{${\bf{K(A)\to{}_{ -\varepsilon} U^{\hbox{{\small{\'el}}}}(A)}}$}}$ qui est
induite par le morphisme d'anneaux $A\times A^{op}\to M_2(A)$ 
d\'ecrit pr\'ec\'edemment. Pour compl\'eter
ce diagramme, nous utilisons un \'el\'ement remarquable de ${}_{-1}D_0
^{\operatorname{max}}(\mathbb{Z})$ et effectuons le ``cup-produit"
par cet \'el\'ement pour d\'efinir une application naturelle 
$\sigma:\bf{{}_\varepsilon K\mathcal{Q}^{\hbox{{\small{\'el}}}}(A) \to{}_{-\varepsilon} D^{
\hbox{{\small{\'el}}}}(A)}$. Les d\'etails sont explicit\'es en \cite{K2} \S  2.3-8 (le fait que 1 soit 
\'eventuellement scind\'e dans $A$ n'est pas n\'ecessaire pour cet argument, comme il a \'et\'e 
d\'ej\`a soulign\'e dans \cite{K2}).

\subsection{\cosa . } Nous proc\'edons de mani\`ere sym\'etrique pour construire une
 application en sens inverse $\bf{{}_{-\varepsilon}U^{\hbox{{\small{\'el}}}}(A) \to
{}_\varepsilon V^{\hbox{{\small{\'el}}}}(SA)}$. Elle s'ins\`ere dans le 
diagramme commutatif suivant
$$\footnotesize{\xymatrix{
\ar[r]&\bf{  \Omega {}_{-\varepsilon} K\mathcal{Q}^{\hbox{{\small{\'el}}}}(A) }\ar[r]\ar[d]&\bf{  {}_{-\varepsilon} U
^{\hbox{{\small{\'el}}}}(A) }\ar[r]\ar[d]&\bf{  K(A) }\ar[r]\ar[d]&\bf{  {}_{-\varepsilon} K
\mathcal{Q}^{\hbox{{\small{\'el}}}}(A) }\ar[r]\ar[d]^\theta&\bf{  {}_{-\varepsilon} U^{\hbox{{\small{\'el}}}}(A)}
\ar[d]\\
\ar[r]&\bf{  {}_\varepsilon E^{\hbox{{\small{\'el}}}}(SA) }\ar[r]&\bf{  {}_\varepsilon V^{\hbox{{\small{\'el}}}}(SA) 
}\ar[r]&\bf{K(A) }\ar[r]&\bf{  {}_\varepsilon E^{\hbox{{\small{\'el}}}}(S^2A) }\ar[r]&\bf{  
{}_\varepsilon V^{\hbox{{\small{\'el}}}}(S^2A)}}}$$
La th\'eorie $\bf{{}_{\varepsilon}E^{\hbox{{\small{\'el}}}}(A)}$ est ici la fibre homotopique de l'application compos\'ee
$$\bf{{}_\varepsilon V^{\hbox{{\small{\'el}}}}(A)\to
{}_\varepsilon V^{\hbox{{\small{\'el}}}}(SA\times SA^{op}) = K(A(e)) \to K(A)}$$
Pour compl\'eter le diagramme, nous devons d\'efinir une application
$$\theta:\bf{{}_{-\varepsilon}K\mathcal{Q}^{\hbox{{\small{\'el}}}}(A) \to  {}_\varepsilon E^{\hbox{{\small{\'el}}}}
(S^2A)}$$
L'id\'ee nouvelle par rapport \`a \cite{K2} est d'utiliser maintenant le cup-produit 
de Clauwens (\'ecrit de mani\`ere relative pour la theorie $E$), soit
$${}_{-1}E_{-2}^{\hbox{\tiny{\'el}}}(\mathbb{Z}[s])\times
{}_{-\varepsilon}K\mathcal{Q}_{n}^{\hbox{\tiny{\'el}}}(A) \to {}_{\varepsilon} E_{n- 2}^{\hbox{\tiny{\'el}}}(A)$$
(avec $\overline{s} = -s$). 

Ceci se traduit au niveau des spectres par l'application $\theta$ . L'\'el\'ement 
de $_{-1}E_{-2}^{\hbox{{\small{\'el}}}}(\mathbb{Z}[s]) =
{}_{-1}K\mathcal{Q}_{-2}
^{\hbox{{\small{\'el}}}}(\mathbb{Z}[s]) ={}_{-1}K\mathcal{Q}_{-2}
^{\operatorname{min}}(\mathbb{Z}[s])$ avec lequel est effectu\'e le cup-produit 
est \'ecrit de mani\`ere
explicite dans \cite{K1} p. 243 par une matrice \`a 30 termes avec un l\'eger changement
de notations (remplacer la lettre ${\lambda}$ par $s$). Nous devons ensuite plonger l'alg\`ebre des polyn\^omes
laurentiens en les deux variables $z$ et $t$ dans la double suspension de $\mathbb{Z}[s]$.

Pour terminer la d\'emonstration du th\'eor\`eme 5.2, nous devons 
montrer que les deux compositions
$$\xymatrix{
\bf{{}_\varepsilon V^{\hbox{{\small{\'el}}}}(SA)}\ar[r]&\bf{{}_{-\varepsilon} U^{\hbox{{\small{\'el}}}}(A)}\ar[r]&
\bf{{}_\varepsilon V^{\hbox{{\small{\'el}}}}(SA)}}\hbox{ et }
\xymatrix{\bf{{}_{-\varepsilon} U^{\hbox{{\small{\'el}}}}(A)}\ar[r]&\bf{{}_\varepsilon V^{\hbox{{\small{\'el}}}}(SA)}
\ar[r]&\bf{{}_{-\varepsilon} U^{\hbox{{\small{\'el}}}}(A)}}$$
sont des \'equivalences d'homotopie. Nous nous r\'ef\'erons de nouveau \`a 
\cite{K2} p 273-277 pour le
d\'etail des arguments. Le point essentiel est l'associativit\'e partielle 
du cup-produit \'etabli en 4.7 qui remplace l'associativit\'e usuelle utilis\'ee 
en \cite{K2}. En effet, de cette associativit\'e partielle,
on d\'eduit des diagrammes commutatifs

$$\xymatrix{
{}_{1}K\mathcal{Q}_0
^{\hbox{{\small{\'el}}}}(\mathbb{Z}[s])
\times {}_\varepsilon K\mathcal{Q}_n
^{\hbox{{\small{\'el}}}}(A)
 \ar[r]\ar[d]&{}_\varepsilon K\mathcal{Q}_n
^{\hbox{{\small{\'el}}}}(A)\ar[d]\\
 {}_{-1}D_0^{\hbox{{\small{\'el}}}} (\mathbb{Z}[s])\times {}_\varepsilon K\mathcal{Q}_n
^{\hbox{{\small{\'el}}}}(A)\ar[r]\ar[d]&
{}_{-\varepsilon} D_n
^{\hbox{{\small{\'el}}}}(A)\ar[d]\\
{}_1K\mathcal{Q}_0
^{\hbox{{\small{\'el}}}}(\mathbb{Z}[s])
\times {}_\varepsilon K\mathcal{Q}_n
^{\hbox{{\small{\'el}}}}(A) \ar[r]&  {}_\varepsilon K\mathcal{Q}_n
^{\hbox{{\small{\'el}}}}(A)
}$$
Ce raisonnement montre que la composition
$$\xymatrix{\bf{ 
{}_\varepsilon V^{\hbox{{\small{\'el}}}}(SA) }\ar[r]&\bf{{}_{-\varepsilon} 
U^{\hbox{{\small{\'el}}}}(A) }\ar[r]&
\bf{{}_\varepsilon V^{\hbox{{\small{\'el}}}}(SA)}}$$
est une \'equivalence d'homotopie. On d\'emontre de m\^eme la commutativit\'e du diagramme

$$\xymatrix{
{}_{-1}D_0
^{\hbox{{\small{\'el}}}} (\mathbb{Z}[s])\times
{}_\varepsilon K\mathcal{Q}_n^{\hbox{{\small{\'el}}}}(A)
 \ar[r]\ar[d]& {}_{-\varepsilon} D_n
^{\hbox{{\small{\'el}}}}(A)\ar[d]\\
{}_1K\mathcal{Q}_0^{\hbox{{\small{\'el}}}}(\mathbb{Z}[s])\times
{}_\varepsilon K\mathcal{Q}_n
^{\hbox{{\small{\'el}}}}(A)
\ar[r]\ar[d]&  {}_\varepsilon K\mathcal{Q}_{n}(A)\ar[d]\\
{}_{-1}D_0
^{\hbox{{\small{\'el}}}} (\mathbb{Z}[s])\times {}_\varepsilon K\mathcal{Q}_n
^{\hbox{{\small{\'el}}}}(A) \ar[r]&  - {}_\varepsilon D_n
^{\hbox{{\small{\'el}}}}(A)
}$$
ce qui montre que la composition en sens inverse
$$\xymatrix{\bf{
{}_{-\varepsilon} U^{\hbox{{\small{\'el}}}}(A)}\ar[r]&\bf{{}_\varepsilon V^{\hbox{{\small{\'el}}}}(SA)}\ar[r]&\bf{{}_
{-\varepsilon} U^{\hbox{{\small{\'el}}}}(A)}
}$$
est aussi une \'equivalence d'homotopie.
\subsection{Remarque.} Si nous nous int\'eressons uniquement aux ``groupes de Witt \'etendus"
$${}_\varepsilon W_n
^{\hbox{{\small{\'el}}}} (A) = Coker(K_{n}(A) \to  {}_\varepsilon K\mathcal{Q}_{n}
^{\hbox{{\small{\'el}}}}(A) )$$
les arguments pr\'ec\'edents se simplifient consid\'erablement (avec un r\'esultat 
moins fort cependant; \`a comparer avec 5.9 et 6.6). Le cup-produit par les \'el\'ements $u_2 \in {}_{-1}W_2
^{\operatorname{max}}(\mathbb{Z})$ et $u_{-2}\in  {}_{-1}W_{-2}
^{\hbox{{\small{\'el}}}}(\mathbb{Z}[s])$, associ\'es aux \'el\'ements construits en 5.5 et 5.6, 
d\'efinissent des homomorphismes
$${}_\varepsilon W_n
^{\hbox{{\small{\'el}}}} (A) \to  {}_{-\varepsilon} W_{n +2}
^{\hbox{{\small{\'el}}}} (A)\hbox{ et }{}_{-\varepsilon} W_{n +2}
^{\hbox{{\small{\'el}}}} (A) \to  {}_\varepsilon W_n
^{\hbox{{\small{\'el}}}} (A)$$
dont la composition (\`a isomorphisme pr\`es) est la multiplication par 4 (en utilisant des
arguments de $K-$th\'eorie topologique : cf. \cite{K1}, p. 251). Notons que ${}_\varepsilon W_n
^{\hbox{{\small{\'el}}}} (A)$ est isomorphe \`a ${}_\varepsilon W_n
^{\operatorname{min}}(A)$ si $n \leq  0$ et \`a ${}_\varepsilon W_n
^{\operatorname{max}}(A)$ si $n < 0$. Le groupe de Witt ``stabilis\'e" que nous
d\'efinirons dans le \S  6 utilisera de mani\`ere essentielle le deuxi\`eme cup-produit.

\subsection{} Comme il a \'et\'e explicit\'e en \cite{K2} p. 278, le 
th\'eor\`eme 5.2 implique une suite exacte \`a 12
termes dont les termes sont d\'efinis ainsi. Le ``cogroupe de Witt" 
${}_\varepsilon\overline{W_n}^{\hbox{{\small{\'el}}}}(A)$ est le noyau de la
fl\`eche oubli
$${}_\varepsilon K\mathcal{Q}_{n}
^{\hbox{{\small{\'el}}}}(A) ) \to K_{n}(A)$$
Nous d\'efinissons le groupe $k_n(A)$ (resp. $\overline{k}_n (A)$ ) 
comme le groupe de cohomologie de Tate pair (resp. impair) de $\mathbb{Z}/2$ 
op\'erant sur $K_{n}(A)$.

{\theorem{}{Avec les d\'efinitions pr\'ec\'edentes, nous avons une suite exacte \`a 12 termes 
o\`{u}, pour simplifier, nous \'ecrivons $F$ pour $F(A)$ en g\'en\'eral, $F$ \'etant l'un des foncteurs
$W^{\hbox{{\small{\'el}}}}$, $\overline{W}^{\hbox{{\small{\'el}}}}$ , $k^{\hbox{{\small{\'el}}}}$ ou 
$\overline{k}^{\hbox{{\small{\'el}}}}$
$$
\begin{array}{c}
\xymatrix{\dots\ar[r]& k_{n+1}
\ar[r]& {}_{-\varepsilon} W_{n +2}
^{\hbox{{\small{\'el}}}} \ar[r]& {}_\varepsilon\overline{W}_n
^{\hbox{{\small{\'el}}}} \ar[r]& \overline{k}_{n+1}
^{\hbox{{\small{\'el}}}} \ar[r]& {}_{-\varepsilon}\overline{W}_{n +1}
^{\hbox{{\small{\'el}}}} \ar[r]& {}_{-\varepsilon} W_{n +1}
^{\hbox{{\small{\'el}}}}}\\
\xymatrix{
\ar[r]& k_{n+1}
\ar[r]& {}_\varepsilon W_{n+ 2}
^{\hbox{{\small{\'el}}}} \ar[r]& {}_{-\varepsilon} \overline{W}_n
^{\hbox{{\small{\'el}}}} \ar[r]& \overline{k}_{n+1}
\ar[r]& {}_\varepsilon \overline{W}_{n+1}
^{\hbox{{\small{\'el}}}} \ar[r]&
{}_\varepsilon W_{n+1}
^{\hbox{{\small{\'el}}}} \ar[r]& k_{n+1}
 \dots} 
\end{array}$$}}

{\theorem{}{ Supposons que 1 soit scind\'e dans $A$ (par exemple que 2 soit inversible).
Les homomorphismes naturels
$${}_\varepsilon W_n
^{\hbox{{\small{\'el}}}} (A) \to  {}_\varepsilon{W}_n(A)\hbox{ et }{}_\varepsilon\overline{W}_n
^{\hbox{{\small{\'el}}}} (A) \to  {}_\varepsilon \overline{W}_n (A)$$
sont alors des isomorphismes.}}
\demo En raisonnant par r\'esurrence sur $n$, c'est une cons\'equence imm\'ediate de 5.9
et du th\'eor\`eme 4.3 de \cite{K2} (voir aussi la remarque suivante).

\subsection{Remarque.} Si 1 est scind\'e dans $A$, nous avons un diagramme commutatif de spectres
$$\xymatrix{
\bf{{}_\varepsilon V^{\hbox{{\small{\'el}}}}(A)}
\ar@{}[r]|-{\hbox{\large{$\approx$}}}&
\bf{\Omega {}_{-\varepsilon} U ^{\hbox{{\small{\'el}}}}(A)}\\
\bf{{}_\varepsilon V(A)}\ar[u]\ar@{}[r]|-{\hbox{\large{$\approx$}}}&\bf{{\Omega}_{-\varepsilon} U (A)}\ar[u]
}$$
o\`{u} les fl\`eches verticales sont des monomorphismes scind\'es. On voit ainsi que le th\'eor\`eme 5.2
implique le th\'eor\`eme fondamental de \cite{K2} p. 260. Nous profitons de cette occasion pour
combler une lacune dans sa d\'emonstration : elle supposait implicitement que 
${}_{1}W(\mathbb{Z}[s])\approx\mathbb{Z}$, un r\'esultat d\^{u} aussi \`a Clauwens (\cite{C} p. 47).

\subsection{} Nous allons conclure ce paragraphe par un calcul explicite de groupes de Witt dans des situations qui ne sont pas envisag\'ees en \cite{K2}. Nous remarquons d'abord que par la m\^eme m\'ethode, nous pouvons d\'efinir en bas degr\'es des morphismes de p\'eriodicit\'e
$$\xymatrix{{}_\varepsilon U^{\hbox{\tiny{min}}}(A)\ar[r]&
{}_{-\varepsilon}V^{\hbox{\tiny{min}}}(SA)}\hbox{ et }\xymatrix{{}_{-\varepsilon}V^{\hbox{\tiny{min}}}(SA)\ar[r]&{}_{\varepsilon}U^{\hbox{\tiny{min}}}(A)}$$
inverses l'un de l'autre \`a isomorphisme pr\`es, en sorte que le diagramme suivant commute
$$\xymatrix{
{}_\varepsilon U^{\hbox{\tiny{\'el}}}(A)\ar[d]\ar[r]&
{}_{-\varepsilon} V^{\hbox{\tiny{\'el}}}(SA)\ar[d]\ar[r]&
{}_{\varepsilon}U^{\hbox{\tiny{\'el}}}(A)\ar[d]\\
{}_\varepsilon U^{\hbox{\tiny{min}}}(A)\ar[r]&
{}_{-\varepsilon} V^{\hbox{\tiny{min}}}(SA)\ar[r]&
{}_{\varepsilon}U^{\hbox{\tiny{min}}}(A)}$$
En effet, la sophistication des fibr\'es plats n'est pas n\'ecessaire dans cette situation. Par ailleurs, puisque 
${}_\varepsilon K\mathcal{Q}^{\hbox{\tiny{\'el}}}(B) $ est isomorphe \`a ${}_\varepsilon K\mathcal{Q}^{\hbox{\tiny{min}}}(B)$ pour tout anneau $B$, on d\'eduit du diagramme pr\'ec\'edent un isomorphisme ${}_\varepsilon U^{\hbox{\tiny{\'el}}}(A)\stackrel{\approx}{\to}
{}_{-\varepsilon} U^{\hbox{\tiny{min}}}(A)$. Nous avons enfin le diagramme commutatif suivant de suites exactes
$$\xymatrix{
0\ar[r]&{}_\varepsilon W_1^{\hbox{\tiny{\'el}}}(A)\ar[d]\ar[r]&
{}_{\varepsilon} U^{\hbox{\tiny{\'el}}}(A)\ar[d]\ar[r]&
K(A)\ar[d]\ar[r]&{}_{\varepsilon}K\mathcal{Q}^{\hbox{\tiny{\'el}}}(A)\ar[d]\\
0\ar[r]&{}_\varepsilon W_1^{\hbox{\tiny{min}}}(A)\ar[r]&
{}_{\varepsilon} U^{\hbox{\tiny{min}}}(A)\ar[r]&
K(A)\ar[r]&{}_{\varepsilon}K\mathcal{Q}^{\hbox{\tiny{min}}}(A)
}$$
Puisque les trois fl\`eches de droite verticales sont des isomorphismes, nous en d\'eduisons le th\'eor\`eme suivant
{\theorem{}{L'homomorphisme naturel
$${}_\varepsilon W_1^{\hbox{\tiny{\'el}}}(A)\to {}_\varepsilon W_1^{\hbox{\tiny{min}}}(A)$$
}}
est un isomorphisme.

\subsection{Exemple} Soit $A=\mathbb{F}_q$ un corps fini de caract\'eristique 2. D'apr\`es Quillen, les groupes $K_n(\mathbb{F}_q)$ sont des groupes finis d'ordre impair \`a l'exception de $K_0(\mathbb{F}_q)=\mathbb{Z}$. On a $W_0(\mathbb{F}_q)=\mathbb{Z}/2$, isomorphisme d\'efini par l'invariant de Arf et $W_1(\mathbb{F}_q)=\mathbb{Z}/2$, isomorphisme d\'efini par l'invariant de Dickson. Ici les groupes de Witt sont ceux calcul\'es avec la forme param\`etre min (c'est-\`a-dire ceux associ\'es \`a des formes \underline{quadratiques}).\\
Par ailleurs, la suite exacte des 12 (th\'eor\`eme 5.11) se r\'eduit en fait \`a une suite \`a 6 termes, car $\varepsilon=1=-1$. Si on utilise le th\'eor\`eme pr\'ec\'edent, on en d\'eduit que les gorupes de Witt \'elargis $W^{\hbox{\tiny{\'el}}}_n(\mathbb{F}_q)$ sont \'egaux \`a $\mathbb{Z}/2$ pour tout $n\in\mathbb{Z}$.

\section{Les groupes de Witt stabilis\'es}
\remark
Ce paragraphe est une extension aux anneaux quelconques des id\'ees d\'evelopp\'ees dans
une Note aux Comptes Rendus \cite{K4}. Une autre extension aux sch\'emas est d\'ecrite dans \cite{S}.

\subsection{} Nous nous pla\c{c}ons dans la cat\'egorie des anneaux discrets $A$ 
avec involution $a {\mapsto}\overline{a}$ (nous ne supposons pas la commutativit\'e ni 
l'existence d'un \'el\'ement unit\'e). Les groupes de Witt stabilis\'es ${}_{\varepsilon}\mathcal{W}_{n}(A)$, 
avec ${\varepsilon} = {\pm} 1$ et $n \in \mathbb{Z}$, que nous d\'efinirons plus loin, v\'erifient les
propri\'et\'es suivantes
\begin{itemize}
\item[\textbf{1)}]\textbf{Exactitude}. Pour toute suite exacte d'anneaux discrets 
avec involution
$$\xymatrix{0\ar[r]&A'\ar[r]&A\ar[r]&A''\ar[r]&0}$$
nous avons une suite exacte naturelle des groupes $\mathcal{W}$
$$\small{\xymatrix{
\ar[r]&  {}_\varepsilon \mathcal{W}_{n+1}(A) \ar[r]&  {}_\varepsilon \mathcal{W}_{n+1}(A") \ar[r]&  {}_\varepsilon 
\mathcal{W}_n(A') \ar[r]&  {}_\varepsilon \mathcal{W}_n(A) \ar[r]&
{}_\varepsilon \mathcal{W}_n(A'') \ar[r]& }}$$
\item[\textbf{2)}]\textbf{Periodicit\'e}. Nous avons un isomorphisme naturel
$${}_\varepsilon  \mathcal{W}_n(A) \cong   {}_{-\varepsilon}  \mathcal{W}_{n+2}(A)$$
et par cons\'equent une p\'eriodicit\'e 4 par rapport \`a l'indice $n$.
\item[\textbf{3)}]\textbf{Invariance par extension nilpotente}. Si $I$ est un 
id\'eal nilpotent dans $A$, la projection $A \to  A/I$ induit un isomorphisme
$${}_\varepsilon  \mathcal{W}_n(A) \cong   {}_\varepsilon  \mathcal{W}_n(A/I)$$
En d'autres termes ${}_\varepsilon  \mathcal{W}_n(I) = 0$ pour un anneau nilpotent.
\item[\textbf{4)}]\textbf{Invariance homotopique}. Si 1 est scind\'e dans $A$ 
(en particulier si 2 est inversible) l'extension polynomiale $A \to  A[t]$ 
(o\`u $\overline{t} = t$ ) induit un isomorphisme
$${}_\varepsilon  \mathcal{W}_n(A) \cong   {}_\varepsilon  \mathcal{W}_n(A[t])$$
\item[\textbf{5)}]\textbf{Normalisation}. Si $A$ est unitaire, il existe un 
homomorphisme naturel
$$\Theta:{}_\varepsilon W_n(A) \to  {}_\varepsilon  \mathcal{W}_n(A)$$
o\`u ${}_\varepsilon W_n(A)$ est le groupe de Witt classique \cite{K1} construit avec les formes quadratiques. 
Celui-ci induit un isomorphisme
$${}_\varepsilon W_n(A) \otimes_\mathbb{Z}\mathbb{Z}'\cong{}_\varepsilon
\mathcal{W}_n(A) \otimes_\mathbb{Z}\mathbb{Z}'$$
o\`{u} $\mathbb{Z}' =\mathbb{Z}[1/2]$.\\
Si $A$ est noeth\'erien r\'egulier, l'homomorphisme ${\Theta}$  est un isomorphisme lorsque $n \leq 0$. Si on
suppose en outre que 2 est inversible dans $A$, les ${}_1\mathcal{W}_n(A)$, $n \mod 4$, sont les groupes 
de Witt triangul\'es de Balmer [Ba].
\end{itemize}

\subsection{} Pour d\'emontrer l'existence d'une telle th\'eorie, 
nous allons essentiellement utiliser les r\'esultats du paragraphe 
pr\'ec\'edent sur la p\'eriodicit\'e en $K$-th\'eorie hermitienne. 
Rappelons que dans \cite{K2} p. 243 nous avons d\'efini un \'el\'ement 
remarquable $u_{-2}$ dans 
$${}_{-1}K\mathcal{Q}_{-2}(\mathbb{Z}[s])=
{}_{-1}K\mathcal{Q}_{-2}^{\hbox{{\small{\'el}}}} (\mathbb{Z}[s])$$ 
d\'efini par une matrice antisym\'etrique ayant 30 \'el\'ements 
et \`a coefficients dans l'anneau des
polyn\^omes laurentiens \`a deux variables $\mathbb{Z}[s][t, u, t^{-1}, u^{-1}]$. 
Cet \'el\'ement nous a d\'ej\`a servi dans le \S 5 pour d\'efinir la fl\`eche ${}_{-\varepsilon} U_{n+1}
^{\hbox{{\small{\'el}}}} (A) \to  {}_\varepsilon V_n
^{\hbox{{\small{\'el}}}} (A)$.

\subsection{} Dans le paragraphe 4, nous avons d\'efini pour tout anneau unitaire $A$ un cup-produit
$${}_{-1}K\mathcal{Q}_{-2}
^{\hbox{{\small{\'el}}}} (\mathbb{Z}[s])\times {}_\varepsilon  K\mathcal{Q}_{n}
^{\hbox{{\small{\'el}}}}(A) \to  {}_{-\varepsilon}  K\mathcal{Q}_{n- 2}
^{\hbox{{\small{\'el}}}} (A)$$
Puisque nous sommes seulement int\'eress\'es aux valeurs de $n$ qui 
sont $\leq  0$, nous pouvons remplacer les groupes 
$K\mathcal{Q}_{n}^{\hbox{{\small{\'el}}}}$ par 
$K\mathcal{Q}_{n}^{\operatorname{min}}$(et m\^eme 
$K\mathcal{Q}_{n}^{\operatorname{max}}$ pour $n < 0$), que nous noterons
simplement $K\mathcal{Q}_{n}$. En outre, l'homomorphisme de p\'eriodicit\'e 
(d\'efini par le cup-produit avec $u_{-2}$)
$$\beta:{}_\varepsilon K\mathcal{Q}_{n}(A)\to
{}_{-\varepsilon} K\mathcal{Q}_{n-2}(A)$$
compos\'e \`a gauche par la fl\`eche oubli 
${}_{-\varepsilon} K\mathcal{Q}_{n-2}(A) \to  K_{n-2}(A)$ ou 
compos\'e \`a droite par la
fl\`eche hyperbolique $K_{n}(A) \to  {}_\varepsilon K\mathcal{Q}_{n}(A)$ 
est r\'eduit \`a 0 (car la $K$-th\'eorie de la suspension d'un
anneau noeth\'erien r\'egulier est triviale). Par cons\'equent, la limite inductive du syst\`eme de
groupes de $K-$th\'eorie hermitienne
$$\xymatrix{
{}_\varepsilon K\mathcal{Q}_{n}(A) \ar[r]&  {}_{-\varepsilon} K\mathcal{Q}_{n-2}(A) \ar[r]&  {}_\varepsilon K\mathcal
{Q}_{n-4}(A) \ar[r]&  {}_{-\varepsilon} K\mathcal{Q}_{n-6}(A) \ar[r]&\dots}$$
est aussi la limite inductive du syst\`eme de groupes de Witt associ\'es
$$\xymatrix{
{}_\varepsilon W_n(A) \ar[r]&  {}_{-\varepsilon} W_{n-2}(A) \ar[r]&  {}_\varepsilon W_{n-4}(A) \ar[r]&  {}_{-
\varepsilon} W_{n-6}(A) \ar[r]&\dots}$$
Cette limite est par d\'efinition le groupe de Witt stabilis\'e 
${}_\varepsilon\mathcal{W}_n(A)$ que nous souhaitions d\'efinir.
Notons que gr\^ace \`a l'excision en $K$-th\'eorie et en $K$-th\'eorie 
hermitienne en degr\'es $\leq  0$, nous pouvons \'etendre cette d\'efinition aux 
anneaux non n\'ecessairement unitaires en
d\'efinissant ${}_\varepsilon K\mathcal{Q}_{n}(A)$ comme le noyau de 
${}_\varepsilon K\mathcal{Q}_{n}(A^+)\to{}_\varepsilon K\mathcal{Q}_{n}(\mathbb{Z})$
, o\`u $A^+$ est l'anneau $A$ (consid\'er\'e comme une $\mathbb{Z}$-alg\`ebre) 
apr\`es addition d'un \'el\'ement unit\'e. La d\'efinition de ${}_\varepsilon\mathcal{W}_n(A)$
pour $A$ non unitaire est tout \`a fait analogue. De ces consid\'erations et 
de l'excision pour les groupes $K\mathcal{Q}_{n}$ si $n \leq  0$, nous d\'eduisons 
la premi\`ere propri\'et\'e des groupes de Witt stabilis\'es:

{\theorem{} {A toute suite exacte d'anneaux discrets avec involution
$$\xymatrix{0 \ar[r]&  A' \ar[r]&  A \ar[r]&  A" \ar[r]&  0}$$
nous pouvons associer naturellement une suite exacte des groupes de Witt stabilis\'es
$$\footnotesize{\xymatrix{
\dots\ar[r]&  {}_\varepsilon  \mathcal{W}_n(A') \ar[r]&  {}_\varepsilon  \mathcal{W}_n(A) \ar[r]&
{}_\varepsilon  \mathcal{W}_n(A'') \ar[r]& {}_\varepsilon \mathcal{W}_{n-1}(A') \ar[r]&{}_\varepsilon \mathcal{W}_{n-1}(A) \ar[r]& \dots}}$$}}

\subsection{} L'isomorphisme ${}_\varepsilon\mathcal{W}_n(A)\cong{}_{-\varepsilon}
\mathcal{W}_{n+2}(A)$ et la p\'eriodicit\'e 4 se d\'eduisent imm\'ediatement des
d\'efinitions.

{\theorem{(normalisation)}{Soit $A$ un anneau noeth\'erien r\'egulier unitaire. Alors le
groupe de Witt stabilis\'e ${}_1\mathcal{W}_0(A)$ (resp.${}_{-1}W_0(A)$) co\"{\i}ncide 
avec le groupe de Witt classique des
formes quadratiques (resp. (-1)-quadratiques). En outre, pour \underline{tout} anneau $A$ unitaire, les
homomorphismes canoniques 
$$\xymatrix{{}_\varepsilon W_n
^{\hbox{{\small{\'el}}}} (A) \ar[r]&  {}_\varepsilon W_n
^{\operatorname{min}}(A) \ar[r]&  {}_\varepsilon W_n
^{\operatorname{max}}(A) \ar[r]&  {}_\varepsilon\mathcal{W}_n(A)}$$ 
induisent des isomorphismes en tensorisant par $\mathbb{Z}' =\mathbb{Z}[1/2]$.}}
\demo Puisque les groupes de $K$-th\'eorie n\'egative de $A$ sont triviaux si $A$ est
noeth\'erien r\'egulier, la suite exacte \`a 12 termes d\'ecrite en 5.9 montre que les fl\`eches de la suite
$$\xymatrix{{}_\varepsilon W_0(A) \ar[r]&  {}_{-\varepsilon} W_{-2}(A) \ar[r]&  {}_\varepsilon W_{-4}(A) \ar[r]&  {}_
{-\varepsilon} W_{-6}(A) \ar[r]&  ...}$$
sont des isomorphismes. Par exemple, si $\varepsilon=1$ et si $A$ est le corps \`a 2 \'el\'ements, nous trouvons
le groupe $\mathbb{Z}/2$ (qui est d\'etect\'e par l'invariant de Arf).\\
Par ailleurs si $A$ est un anneau quelconque, en utilisant la localisation en $K$-th\'eorie hermitienne, nous avons 
construit en \cite{K1}
deux \'el\'ements dans ${}_{-1}W_2
^{\operatorname{max}}(\mathbb{Z})$ et ${}_{-1}W_{-2}
^{\operatorname{max}}(\mathbb{Z})$ dont le cup-produit dans ${}_1W^{\operatorname{max}}(\mathbb{Z})$ est une
puissance de 2. Les premiers isomorphismes se d\'emontrent en se ramenant par p\'eriodicit\'e aux
degr\'es n\'egatifs. Le dernier isomorphisme r\'esulte de la suite exacte \`a 12 termes d\'emontr\'ee en
5.9.

{\theorem {(invariance par extension nilpotente)} {Si $I$ est un id\'eal nilpotent dans $A$, la
projection $A \to  A/I$ induit un isomorphisme
$${}_\varepsilon\mathcal{W}_n(A)\cong{}_\varepsilon\mathcal{W}_n(A/I)$$
Par cons\'equent, ${}_\varepsilon\mathcal{W}_n(I)=0$ pour tout id\'eal nilpotent $I$.}}
\demo Sans restreindre la g\'en\'eralit\'e, nous pouvons supposer que $A$ est unitaire. Dans
ce cas, il est bien connu que tout module projectif de type fini sur $A/I$ provient d'un module
projectif $E$ sur $A$ par extension des scalaires et qu'il est donc du type $E/I$. Par cons\'equent, la
forme $\varepsilon$-hermitienne sur $A/I$ est donn\'ee par un isomorphisme
$$\varphi  : E/I \to  (E/I)^*$$
Puisque $\varphi$ est paire, nous pouvons l'\'ecrire sour la forme $\varphi_0 + \varepsilon{}^ t\varphi_0$. 
Soit $\widetilde{\varphi}_0$ un homomorphisme $E \to E^*$ tel que $\widetilde{\varphi}_0=\varphi_0 \mod I$. 
Alors $\varphi  =\widetilde{\varphi_0}+\varepsilon{}^t\widetilde{\varphi_0}$ est une forme $\varepsilon$-
hermitienne non d\'eg\'en\'er\'ee $E \to E^*$ qui est un relev\'e de $\varphi$. Ceci montre que le morphisme
${}_\varepsilon K\mathcal{Q}_{0}(A) \to  {}_\varepsilon K\mathcal{Q}_{0}(A/I)$ 
est surjectif pour tout id\'eal nilpotent $I$ (aussi bien pour $K\mathcal{Q}^{\operatorname{max}}$ que
pour $K\mathcal{Q}^{\operatorname{min}}$). Il en est donc de m\^eme de
$${}_\varepsilon K\mathcal{Q}_{n}(A) \to  {}_\varepsilon K\mathcal{Q}_{n}(A/I)$$
pour $n \leq  0$ en consid\'erant des suspensions it\'er\'ees (\'el, max et min co\"{\i}ncident en degr\'es $n < 0$;
cf. 3.2). La surjectivit\'e de l'homomorphisme ${}_\varepsilon\mathcal{W}_n(A)
\to{}_\varepsilon\mathcal{W}_n(A/I)$ en r\'esulte.\\
L'injectivit\'e du morphisme ${}_\varepsilon\mathcal{W}_n(A)\to{}_\varepsilon\mathcal{W}_n(A/I)$
 est plus d\'elicate \`a montrer. En raisonnant par r\'ecurrence sur le degr\'e de nilpotence de $I$, 
nous pouvons d'abord supposer que $I^2 = 0$. Par ailleurs, nous savons que tout module muni d'une 
forme hermitienne paire est facteur direct d'un module hyperbolique. C'est donc l'image d'un projecteur 
auto-adjoint $p$, soit $p^2=p$ et $p^*=p$ dans un $H(A^n)$.\\
Enfin, sans restreindre la g\'en\'eralit\'e (puisque nous stabilisons), nous pouvons supposer que $A$
est la suspension $SR$ d'un anneau $R$ et que $I = SJ$ o\`u $J$ est un id\'eal de $R$ tel que $J^2 = 0$.
La d\'emonstration de l'injectivit\'e se r\'esume alors \`a la solution du probl\`eme suivant : nous
consid\'erons deux projecteurs auto-adjoints $p_0$ et $p_1$ dans un module hyperbolique sur $A = SR$
tels que leurs images $\mod\ I$, soient $\overline{p}_0$ et $\overline{p}_1$ sont conjugu\'ees. Puisque 
${}_\varepsilon K\mathcal{Q}_{1}(SR) \cong   {}_\varepsilon K\mathcal{Q}_{0}(R)$ en
g\'en\'eral et que le morphisme ${}_\varepsilon K\mathcal{Q}_{0}(R)\to{}_\varepsilon K\mathcal{Q}_{0}(R/J)$
 est surjectif comme nous l'avons vu pr\'ec\'edemment, nous pouvons supposer sans restreindre 
la g\'en\'eralit\'e\footnote{La surjectivit\'e de l'homomorphisme ${}_\varepsilon K\mathcal{Q}_{1}
^{\operatorname{min}}(\Lambda ) \to  {}_\varepsilon K\mathcal{Q}_{1}
^{\operatorname{min}}(\Lambda  / I )$ implique la surjectivit\'e de
l'homomorphisme ${}_\varepsilon O^{\operatorname{min}}( \Lambda ) \to  
{}_\varepsilon O^{\operatorname{min}}(\Lambda /I)$.} que $\overline{p}_0=\overline{p}_1$ ou encore
$p_1=p_0+\sigma$, o\`u $\sigma$ appartient \`a $I$. De l'identit\'e $(p_1)^2 = p_1$ et de l'\'egalit\'e $I^2 = 0$, nous 
d\'eduisons les relations suivantes :
\begin{eqnarray*}
{\sigma}  = p_0 {\sigma}  + {\sigma}  p_0\\
{\sigma}  p_0 {\sigma}  = 0\\
{\sigma}^2 = {\sigma}^2p_0 = p_0{\sigma}^2
\end{eqnarray*}

Consid\'erons maintenant l'endomorphisme $\alpha  = 1 - p_0 - p_1 + 2 p_0 p_1$. Puisque $\alpha  \equiv  1 \mod I$,
c'est un isomorphisme. Par ailleurs, il v\'erifie la relation $\alpha p_1=p_0\alpha$.
Nous allons maintenant montrer que $\alpha  \alpha^* = 1$. Pour cela, on remarque que $\alpha$ s'\'ecrit aussi
$$\alpha  = 1 - \sigma  +2p_0\sigma$$ 
et, gr\^ace aux identit\'es pr\'ec\'edentes, un calcul direct montre bien que 
$$\alpha \alpha^* = (1 - \sigma  +2p_0\sigma)(1 - \sigma  + 2\sigma p_0) = 1$$
Les projecteurs $p_0$ et $p_1$ sont ainsi conjugu\'es par un automorphisme unitaire et d\'eterminent par
cons\'equent la m\^eme classe de forme hermitienne paire\footnote{D'apr\`es 2.10, il revient au m\^eme de consid\'erer des formes hermitienne paires 
ou des formes quadratiques dans les groupes stabilis\'es.}.

{\theorem {(invariance homotopique)} {Soit $A$ un anneau unitaire tel que 1 soit scind\'e
dans $A$. Il existe donc un \'el\'ement $\lambda$ dans le centre de $A$ tel que 
$\lambda +\overline{\lambda}=1$. L'extension
polynomiale $A \to  A[t]$ (avec $\overline{t}=t$) induit alors un isomorphisme
$${}_\varepsilon\mathcal{W}_n(A) \cong {}_\varepsilon\mathcal{W}_n(A[t])$$}}
\demo Il suffit de d\'emontrer le th\'eor\`eme pour $n = 0$. Celui-ci est d\'ej\`a connu pour 2
inversible dans $A$ (voir \cite{O} pour une preuve simple). Cependant, il existe des anneaux o\`u 2
n'est pas inversible et o\`u 1 est scind\'e, par exemple le corps fini $\mathbb{F}_4$ muni de l'involution non
triviale. Pour traiter ce cas plus g\'en\'eral, nous devons r\'e\'examiner la preuve classique. En fait, le
seul point qui m\'erite une pr\'ecision dans cette preuve est le lemme suivant.

{\lemma{Soit $A$ un anneau avec $\lambda$  dans le centre de $A$ tel
que $1 = \lambda  +\overline{\lambda}$. Soit $E$ un $A$-module muni d'une forme $\varepsilon$-hermitienne 
et soit $\alpha  = 1 + \nu  t$
un \'el\'ement de $GL(E \otimes\mathbb{Z}[t])$ avec $\nu$ nilpotent et auto-adjoint. Alors $\alpha$ 
peut \^etre \'ecrit sous la forme $\gamma (t)^*\gamma (t)$, o\`u $\gamma(t)$ est un polyn\^ome en $t$ 
dans l'anneau engendr\'e par $\lambda$  et $\nu$.}}
\demo Nous allons construire par r\'ecurrence sur $n$ un polyn\^ome de degr\'e au plus $n$
dans l'anneau engendr\'e par $\nu$ et $\lambda$, soit $\gamma_n(t) =1 + a_1t + a_2t^2 +\dots + a_nt^n$, 
tel que $\gamma_n(t)^*\gamma_n(t)\equiv  1 + \nu  t \mod (\nu  t)^{n+1}$. Pour $n = 1$, nous posons 
$\gamma_1(t) = 1 + \lambda  \nu  t$. Si $\gamma_n$ est construit, nous
avons $\gamma_n(t)^*\gamma_n(t) = 1 + \nu  t + b_{n+1} (\nu t)^{n+1} \mod (\nu t)^{n+2}$ avec 
$b_{n+1}=\overline{b}_{n+1}$. Nous posons
alors $\gamma_{n+1}(t)=(1-\lambda b_{n+1}(\nu t)^{n+1})\gamma_n(t)$ pour obtenir l'identit\'e requise
$$\gamma_{n+1}(t)^*\gamma_{n+1}(t) \equiv  1 + \nu  t \mod (\nu t)^{n+2}$$

\subsection{Exemple} Si A est un corps fini de caract\'eristique 2, il est facile de montrer que les
groupes de Witt stabilis\'es $\mathcal{W}_{n}(A)$ sont tous isomorphes \`a $\mathbb{Z}/2$. Ils co\"{i}ncident 
en fait avec les groupes ${W}_n^{\hbox{\tiny{\'el}}}(A)$ en tout degr\'e.

\remark Ces groupes de Witt stabilis\'es ont \'et\'e g\'en\'eralis\'es aux sch\'emas par M.
Schlichting \cite{S}. Dans cette g\'en\'eralit\'e, on doit cependant supposer 2 inversible.

\section[APPENDICE]{Les lemmes de Clauwens}
{\lemma{La forme hermitienne associ\'ee \`a la forme quadratique ${\kappa}$  d\'efinie en 4.3 est non
d\'eg\'en\'er\'ee.}}
\demo Nous suivons les simplications de notation indiqu\'ees en 4.3 en rempla\c{c}ant
notamment ${\delta}$ par ${\phi}$  tel que ${\phi}  + {\phi}^* = 1$. Nous pouvons donc \'ecrire
$${\kappa}  = \sum  {\theta}_n {\otimes} {\phi}^n$$
qu'il est plus suggestif de noter ${\theta}({\phi} )$.
Nous avons alors 
$${\kappa}  + {\kappa}^* =\sum  {\theta}_n {\otimes} {\phi}^n +\sum  ({\theta}_n)^* {\otimes} ({\phi} )^{*n}=
\sum  {\theta}_n {\otimes} {\phi}^n +\sum  ({\theta}_n)^* {\otimes} (1 - {\phi} )^n$$
Par ailleurs, on sait que le polyn\^ome en $s$ d\'efini par $\sum  {\theta}_n {\otimes} s^n +
\sum({\theta}_n)^* {\otimes} (1 - s)^n$
est inversible (c'est la forme hermitienne $H$ associ\'ee \`a ${\theta}$). Il en r\'esulte \'evidemment que ${\kappa}  + {\kappa}^*$ est inversible. On peut aussi l'\'ecrire $H({\phi} )$ avec un abus d'\'ecriture \'evident.
{\lemma {Si on change ${\theta} = \sum  {\theta}_n s^n$ en ${\theta} + Z - Z^*$, les formes quadratiques 
associ\'es ${\kappa}$  et ${\kappa}'$
sont \'equivalentes.}}
\demo La forme quadratique ${\theta} =\sum  {\theta}_n s^n$ est modifi\'ee en
$$\sum  {\theta}_n s^n + \sum  {\sigma}_n s^n - \sum  ({\sigma}_n)^* (1-s)^n$$
Par cons\'equent ${\kappa}$ est modifi\'ee en ${\kappa}  + {\sigma} ({\phi} ) - ({\sigma} ({\phi} ))^*$ 
(remplacer $s$ par ${\phi}$ ).
{\lemma{Modulo l'image de $K\mathcal{Q}(A)$ dans $K\mathcal{Q}^{\hbox{\tiny{\'el}}_0}(A[s])$ 
(et m\^eme d'une forme hyperbolique sur $A$), tout \'el\'ement de ce dernier groupe peut \^etre repr\'esent\'e 
par une forme lin\'eaire en $s$.}}
\demo Soit ${\theta} = \sum_0^N{\theta}_ns^n$ une forme quadratique de degr\'e $N$. L'identit\'e suivante et un
raisonnement par r\'ecurrence sur $N$ montre qu'on peut r\'eduire le degr\'e de ${\theta}$ \`a 0 ou 1
$$\left(\begin{array}{ccc}
1&-s& ({\theta}_N )^*(1-s)^{N-1}\\
0& 1& 0\\
0& 0& 1
\end{array}\right)
\left(
\begin{array}{ccc}
{\theta}& 0& 0\\
0& 0&1\\
0 &0 &0
\end{array}\right)
\left(
\begin{array}{ccc}
1& 0 &0\\
-1+ s& 1& 0\\
{\theta}_Ns^{N-1}& 0& 1
\end{array}\right)=
\left(
\begin{array}{ccc}
{\theta} -{\theta}_N s^N &0& -s\\
{\theta}_Ns^{N-1}& 1 &0\\
0&0& 1
\end{array}\right)$$

Si ${\theta}$ s'\'ecrit ${\theta}_0 + {\theta}_1s$, on peut aussi \'eliminer le terme constant en \'ecrivant que 
${\theta}$ est \'equivalente \`a $$	{\theta}_0 + {\theta}_1s - {\theta}_0(1-s) + ({\theta}_0)^*s = 
({\theta}_1 + {\theta}_0 + ({\theta}_0)^*)s$$
ce qui d\'emontre le lemme.\\ \\
Le lemme pr\'ec\'edent nous montre qu'il suffit de v\'erifier la validit\'e du produit de Clauwens
d\'efini en 4.3 (m\^emes notations), dans le cas o\`{u} ${\theta}$ est une forme lin\'eaire en s, soit 
${\sigma} s$ avec ${\sigma}$	 presque sym\'etrique, i.e. ${\sigma}^* = {\sigma} (1+N)$, avec $N$ nilpotent. 
Il nous faut montrer ensuite que le cup-produit de Clauwens ne d\'epend que de la forme 
quadratique associ\'ee \`a ${\delta}$ (ou
l'endomorphisme ${\phi}$  gr\^{a}ce \`a l'identification de $F$ \`a son dual). Rappelons qu'on a aussi identifi\'e
$E$ \`a son dual par l'isomorphisme ${\theta}_0 + \sum_{n=0}^{\infty}{}^t{\theta}_n$.\\ \\
Si on pose $G = E {\otimes} F$, la transpos\'ee ${}^{t}f$ d'une application $f$ de $G$ dans son dual s'identifie
\'egalement \`a son application adjointe $f^*$ (cf. les remarques faites en 4.3).

{\lemma{Soit ${\phi}$  et ${\zeta}$ deux endomorphismes de $F$ tels que ${\phi}  + {\phi}^* = 1$. Pour tout entier 
$p \geq  0$, il existe alors un isomorphisme $f_p$ de $G$ sur son dual tel que 
$$(f_p)^*({\sigma} {\otimes} {\phi} )f_p = {\sigma}  {\otimes} ({\phi}  + {\zeta} - {\zeta}^*) + Z_p - (Z_p)^* 
\mod ({\sigma} N^{p+1} {\otimes} 1)$$
o\`{u} $N = {\sigma}^{-1} {\sigma}^*-1$ est nilpotent et o\`{u} l'expression mod 
$({\sigma} N^{p+1} {\otimes} 1)$ signifie une somme de morphismes du type 
${\sigma} N^{p+1} {\otimes} {\kappa}_{p+1}+{\sigma} N^{p+2}{\otimes} {\kappa}_{p+2}+\dots$ 
(qui est finie puisque $N$ est nilpotent).}}
\demo Puisque ${\sigma}^* = {\sigma}+{\sigma}N$, on a ${\sigma}^* N^k {\otimes} 1 = 
{\sigma}  N^k {\otimes} 1 \mod ({\sigma}  N^{k+1} {\otimes} 1)$.
On a de m\^eme $N^{*k}{\sigma} N^r = {\sigma} N^{r+k}\mod {\sigma}N^{r+k+1} {\otimes} 1$. 
Nous allons maintenant construire $f_p$ et $Z_p$ par r\'ecurrence sur $p$. 
Pour $p=0$, on pose $f_0 = 1$ et $Z_0 = - {\sigma} {\otimes} {\zeta}$. 
Pour d\'efinir $f_{p+1}$ \`a partir de $f_p$ , on \'ecrit
$$(f_p)^*({\sigma}  {\otimes} {\phi} )f_p - [ {\sigma}  {\otimes} ({\phi}  + {\zeta} - {\zeta}^*) + Z_p - (Z_p)^* ] = - {\sigma}  N^{p+1} {\otimes} {\kappa}_{p+1}\mod ({\sigma}N^{p+1} {\otimes} 1)$$
On pose alors
$U = N_{p+1} {\otimes} {\kappa}_{p+1}$ et $f_{p+1} = f_p + U$ et $Z_{p+1} = Z_p + U^* ({\sigma}  {\otimes} {\phi} )$\\
En travaillant mod $({\sigma} N^{p+2} {\otimes} 1)$, on obtient les identit\'es suivantes

$(f_{p+1})^*({\sigma}  {\otimes} {\phi} )f_{p+1} - [ {\sigma}  {\otimes} ({\phi}  + {\zeta} - {\zeta}^*) + Z_{p+1} - (Z_{p+1})^* ]$\\
$= (f_{p+1})^*({\sigma}  {\otimes} {\phi} )f_{p+1} - (f_p)^*({\sigma}  {\otimes} {\phi} )f_p - (Z_{p+1} - Z_p) + ((Z_{p+1})^* - (Z_p)^*) - {\sigma}  N^{p+1} {\otimes} {\kappa}$ 
$= U^*({\sigma}  {\otimes} {\phi} ) + ({\sigma}  {\otimes} {\phi} )U - U^*({\sigma}  {\otimes} {\phi} )) + ({\sigma}^* {\otimes} {\phi}^*)U - {\sigma} N^{p+1} {\otimes} {\kappa} $
$= {\sigma}  N^{p+1} ({\phi}  + {\phi}^* - 1) {\kappa} $
$= 0 \mod {\sigma}  N^{p+2}{\otimes} 1$\\
Ceci ach\`eve la d\'emonstration du lemme.

${}$\\ 
\flushleft{MAX KAROUBI\\
Universit\'e Denis Diderot / Paris 7\\
UFR de Math\'ematiques. Case 7012\\
175, rue du Chevaleret\\
75205 Paris cedex 13\\
\href{mailto:max.karoubi@gmail.com}{max.karoubi@gmail.com}\\
URL: \href{http://people.math.jussieu.fr/%7Ekaroubi}{http://people.math.jussieu.fr/$\sim$karoubi}}
\end{document}